
\input amssym.tex
\magnification=1095
\parindent=0pt

\def\MaxBar{{\underline{\rm Max}\,}}
\def\Max{{\rm Max\,}}
\def\O{{\cal O}}
\def\Sing{{\rm Sing}}
\def\word{\hbox{\rm w-ord}}

\def\frac#1#2{{#1\over #2}}
\font\twelverm=cmr12
\font\twelvebf=cmbx12
\font\twelveit=cmti12
\font\twelvesl=cmsl12
\font\twelvemus=cmmi12
\def\bigtype{\let\rm=\twelverm \let\bf=\twelvebf \let\it=\twelveit
  \let\sl=\twelvesl \let\mus=\twelvemus
  \baselineskip=14pt minus 1pt
  \rm}
\font\eightrm=cmr8


\centerline{\bigtype\bf On Algorithmic Equiresolution and Stratification
of Hilbert Schemes.}
\medskip
\footnote{}{Mathematics Subject Classification 2000: 14E15, 14D99.}

\centerline{S. Encinas\footnote{*}{\eightrm Partially supported by PB96-0065},
A. Nobile\footnote{**}{\eightrm Partially supported by the program
``Estancia de investigadores en r\'egimen de a\~no s\'abatico en
Espa\~na''}, O. Villamayor$^{*}$.}
\bigskip

\centerline{\bf Abstract}
\medskip

\vbox{\leftskip=2cm\rightskip=2cm\eightrm
Given an algorithm of resolution of singularities satisfying
certain conditions (``good algorithms''), natural notions of
simultaneous algorithmic resolution, or equiresolution,  for
families of embedded schemes (parametrized by a reduced scheme $T$)
are proposed. It is proved that these conditions are equivalent.
Something similar is done for families of sheaves of ideals, here
the goal is algorithmic simultaneous principalization. A
consequence is that given a family of embedded schemes over a
reduced $T$,  this parameter scheme can be naturally expressed  as a
disjoint union of locally closed sets $T_{j}$, such that the
induced family on each part $T_{j}$ is
equisolvable. In particular, this can be applied to the Hilbert
scheme of a smooth projective variety; in fact, our result shows that,
in characteristic zero, the underlying topological space of
any Hilbert scheme parametrizing embedded schemes
can be naturally stratified in equiresolvable families.
}
\bigskip

\bigskip

\centerline{\bf INTRODUCTION}
\medskip

By now, the theory of desingularization of varieties in
characteristic zero seems to be pretty well understood. The first
general result on resolution, valid for any dimension,
is the famous theorem of Hironaka,
published in 1964 [15]. This is a precise result, showing not only
that desingularization is possible, but also that given a variety,
or more generally an excellent scheme $X$ over a field of
characteristic zero it is possible, by means of a finite sequence of
blowing-ups, each one with center a regular
subscheme, to obtain a birational, projective morphism $f: X' \to
X$, with $X'$ regular.
Moreover, one gets ``embedded desingularization''. This means: if
$X$ is a subscheme of a regular scheme $W$, it can be reached a
situation where $X'$ is a subscheme of a regular $W'$, and $f$ is
induced by a projective  morphism $g:W' \to W$, inducing an
isomorphism from
$W'-g^{-1}(S) $ onto $W-S$, where $S$ is the singular set of $X$,
in such a way that $X'$ has normal crossings with the exceptional
divisor of $g$ (see 1.1.e of this article for the precise
definition). A possible ``shortcoming'' of this important result is
the fact that, although the morphism $f$ (or $g$, in the embedded
situation) is a composition of ``nice'' blowing-ups, it is not
specified how to choose the center each time. More recent work shows
how to obtain more constructive proofs. We refer here to [6], [10] and
[28] for different algorithms of desingularization. A more
self-contained presentation of the latter appears in [11].
These algorithms tell us, given a
subvariety $X$ of a regular variety $W$ (over a field of
characteristic zero), how to choose the different centers in order
to get embedded desingularization, by repeatedly blowing-up along the
centers provided by the algorithm.

Let us mention that at least the algorithm treated in [11] can be
``implemented''. In fact, in [8] a number of explicit
examples (involving surfaces in tree-space) are resolved by using that
algorithm, with the aid of computers (see also [9]).
\smallskip

The development in [11] also led to a short and simplified proof of
desingularization avoiding Hironaka's notion of normal flatness (see
[12] and addendum in [11]).
\smallskip

It should be mentioned other recent proofs (short, but
non-constructive) obtained with very different techniques,
namely the use of the theory of moduli of curves. See  [17], [5],
[7], [4] and [3].
\smallskip

Once the problem of resolving, in an explicit (or constructive)
way, a {\it single} algebraic variety  has been settled
it is natural to consider the question of ``classification'' of
varieties according to the resolution of their singularities. A
closely related  question is the study of criteria saying when a
family of varieties
can be {simultaneously}
resolved, in an appropriate sense.  The present article is
concerned with this type of questions.

Simple examples indicate that some restrictions have to be imposed
on the family. It seems that the first systematic, rigorous analysis
of the problem, in the case of families of curves, is due to
Zariski, as part of his program to study equisingularity in
codimension one. This work was expanded by other mathematicians,
specially Teissier (see [27] and the bibliography therein). One
of the main results is, essentially: working with complex analytic
varieties, a family of plane curves can be simultaneously
desingularized if and only if certain numerical invariants of the
members of the family remain constant. See [27], Theorem 5.3.1,
for a precise statement (and the proof). Also in the early
nineteen-seventies Risler investigated the analogous problem in
dimension zero, i.e., that of families of ideals in the plane,
having zero-dimensional support (again in the context of local
complex analytic geometry). There was (by the same time, i.e. the
mid 70's) some interesting work on families of surfaces ([19],
[30],) but apparently during the two following decades little
work has been done in this direction. Interest in this area was
revived by Lipman's article [20]. Perhaps the articles [24],
[25], [23] and [29] should also be mentioned. The first two deal
with Risler's theory on families of zero dimensional ideals, the
third with families of curves, the last one with a problem posed by
Lipman in [20], based on questions of Zariski.

The main difficulty we face when we want to move from a family of
varieties whose members have dimension $d\leq 1$ to the case $d \geq
2$ is the lack of ``canonicity'' in the resolution process. In
fact, to resolve a curve $C$, it is ``obvious'' what to do: blow-up
the singular locus $S$ of $C$, to obtain a curve $C'$, with singular
locus $S'$, now blow up $C'$ with center $S'$, and so on. (What is
not so obvious is to prove that this process actually works, as the
great forerunners of the nineteenth and early twentieth century soon
discovered; they had to create the pertinent machinery to cope with
the problem.) Something similar happens in dimension zero, i.e. for
families of zero-dimensional ideals on regular surfaces. But in
higher dimension such a natural choice for the centers disappears. A
remedy is provided by a desingularization algorithm. Of course, the
price to pay is that now we must fix, beforehand, such an
algorithm. This is an ``arbitrary'' choice, but once this is done,
we have a definite selection of the center to use each time we
blow-up, i.e. we are essentially in the same situation as in
dimension one. This is the point of view followed in the present
article.
\smallskip

So, our basic concern is to simultaneously resolve families of
varieties (or suitable schemes over a field of characteristic zero),
when  a resolution algorithm is fixed. However, it seems that,
rather than directly investigating this question, it is convenient
to study a more general, algebraic one, namely to study  what
happens when we have a family of regular varieties (or suitable
schemes) $W_t$ ($t$ ranging in a suitable parameter space $T$), and
for each $t$ a sheaf of ideals ${\cal I}_t$ (for short, a {\it
family of ideals}.) For a single object, i.e. when $T$ is one point,
the goal is to ``principalize'' the sheaf, i.e., by taking a
sequence of blowing-ups and total transforms of the ideal, reach a
situation where the transformed ideal is locally principal. In the
family case, we try to do this via a global process that induces,
fiber-wise, the principalization of each member of the family.
(Again assuming that we have fixed an algorithm of principalization,
which allows us to talk about ``the'' principalization process of
an ideal.)

A little more precisely, the problem is to give ``reasonable''
notions of {\it simultaneous principalization}, or equiresolution,
of a family of (sheaves of) ideals, and to show that these are
equivalent. Ditto for families of embedded varieties (or suitable
schemes.) It seems convenient to attack this ``simultaneous
principalization problem'' first, for several reasons.

First,  this problem (for ideals) is more general than the
analogous problem for embedded varieties. In fact, if properly
formulated, a solution thereof  easily implies a solution of the
simultaneous  embedded resolution problem for varieties.

Second, traditionally the notion of embedded resolution applies to
a pair reduced subscheme-regular ambient scheme. However, even if we
were primarily interested in classical algebraic varieties, when
considering families we cannot avoid the presence of certain members
(``degenerations'') which will be non-reduced schemes. We'll see
that a good notion of simultaneous principalization of a family of
ideals leads to a reasonable notion of simultaneous resolution of
embedded schemes, even if some members of the family are not
reduced.

Finally, it seems that the  more general problem of
principalization is technically no more difficult that that of
resolution, so we get more general results with no added effort.

So, in this article we prioritize the study of families of ideals.
Let us briefly explain, more precisely, the contents of this paper.
In chapter 1, after introducing some notation and terminology, we
define a notion  of {\it strong principalization algorithm} (for a
sheaf of ideals over a suitable regular scheme $W$) and (given a
reduced subscheme $X$ of a suitable regular scheme $W$) that of {\it
desingularization algorithm}. (The meaning of the term ``suitable''
is explained.) We show how a desingularization algorithm easily
follows if we have a principalization one for ideals. We introduce
notions of {\it family of ideals} and {\it family of embedded
schemes}. To deal with equiresolution, we find it necessary to
introduce more conditions on our algorithm, this leads to the
concept of {\it good principalization algorithm}, which is explained
next. Finally, we indicate an example of a good principalization
algorithm, based on the work of [11]. It is not clear that some
of the properties of a good algorithm (namely, 1.13 (4) and (5)) follows
from the results of [11], the validity of this properties is
proved in 1.18 and chapter~6 respectively.
In chapter~2 we introduce two notions of
simultaneous principalization, or ``equiresolution'', for families
of ideals, once a strong principalization algorithm has been fixed.
We call them condition AE (after ``algorithmic equiresolution'') and
condition $\tau$ respectively. We believe both are, in a sense,
natural: the first essentially says that  we may resolve using
centers ``evenly spread'' over the parameter space, the second
imposes the constancy of certain ``numerical'' invariants associated
to each fiber. Next we state a theorem saying that if our algorithm
is good (and some other mild properties are satisfied) then both
notions are equivalent. The implications for families of embedded
schemes are indicated. In particular, we propose a notion of
equisolvability of a family of embedded scheme which allows some of
the fibers to be non-reduced.

In chapter 3 we prove the theorems announced in the previous
chapter. In chapter 4 we prove, among other things, a theorem that
says, essentially, that given a family of embedded schemes,
parametrized by a reduced scheme $T$, it is possible to naturally
express $T$ as a finite disjoint union of locally closed subschemes
$T_j$ such that the restriction of the family to each $T_j$ is
equisolvable. This follows from more general similar results for
families of ideals and it has an application to the universal families
obtained with the aid of the Hilbert Scheme Theory. In this way we
obtain families of subschemes of an ambient projective scheme,
equisolvable and universal with respect to this property.

As indicated
above, the last two chapters are devoted to complete  the proof that
the algorithm of [11] is a good one. This is probably the most
technical and less self-contained part of this article; although we
tried to recall the basic facts and terminology, it seems
unavoidable, for reasons of space, to assume some familiarity with
[11].

Finally, some informal comments on the results of chapter 4. We
believe that a reasonable desingularization algorithm should be
``algebraic'', in the sense that if it is applied to the total space
of, say, an algebraic family of varieties, it should behave
``nicely'' with respect to the fibers (recall that our algorithm
associates a resolution to each of them.)  The result of chapter 4
on the possibility to find a stratification of the parameter space
into locally closed sets with the properties indicated above seems
to be make this vague notion precise, and to show that this property
is valid for good algorithms (in the sense of chapter 1); in
particular for that of [11]. Probably other algorithms
available in the  literature share this property, although we did
not verify this fact.
\smallskip

In general, we use the notation and terminology of [14].
\smallskip

The second author wishes to thank the Ministry of Education and
Culture of Spain for its support while he was  visiting the
Universidad Aut\'onoma de
Madrid in the Fall of 1998; during this period the present
investigation was started and considerably developed.
\bigskip

{\bf 1. GOOD RESOLUTION ALGORITHMS.}
\medskip

{\bf (1.1)} We  shall work with certain ``resolution algorithms'',
which will be defined on suitable classes of schemes. Before
precisely  introducing  such classes, we need a definition.
\smallskip

(a) A morphism $f: Z' \to Z$ of noetherian schemes is a {\it
localization morphism} if, up to isomorphism, it is of the following
form: there must be a fiber product square of schemes and morphisms

$$\displaylines{
\matrix{{Z'}&{\buildrel f \over \rightarrow}&Z\cr
        {{\scriptstyle p'} \downarrow ~}& &{{\scriptstyle
p}\downarrow ~}\cr
          {T'}&{\rightarrow }&T\cr}\cr}$$
where  $T'= {\rm Spec}\,{\cal O}_{T,t}$, for a suitable point $t \in
T$ and $T' \to T$ is the canonical morphism. Note that $f$ is flat
and induces a homeomorphism of the underlying topological space of
$Z'$ with its image (regarded as a topological subspace of $Z$, with
the relative topology).
\smallskip

(b) An {\it allowable collection} is a class ${\cal S}$
such that each $Z \in {\cal S}$ is a
noetherian, regular, excellent scheme over some field $ k $ of
characteristic
$0$ (the
base field
is allowed to vary), satisfying the following conditions:
\smallskip

(1) If $Z \in {\cal S}$, $Z'$ is a regular noetherian scheme  and
there is a localization  morphism $Z' \to Z$, then $Z' \in {\cal
S}$.

(2) If $Z' \to Z$ is a morphism of finite type of regular schemes
and $Z \in {\cal S}$, then
$Z' \in {\cal S}$.

(3) If  $Z \in {\cal S}$ is a scheme over a field $k$, $K$ is an
extension of $k$ and $Z_K$ is the $K$-scheme obtained from
$Z$ by base change, then $Z_K$ is in ${\cal S}$.
\smallskip

(c) An {\it S-pair} is a pair $(W, E)$, where $W \in {\cal S}$ is
regular, equidimensional and $E=\{H_1, \ldots, H_m \}$ is a
collection of regular hypersurfaces of $W$ having normal crossings
(see part (e).)

By a {\it hypersurface} of $W$ we mean a reduced, purely
1-codimensional subscheme of $W$.

A subscheme $C \subset W$ is a {\it permissible center} for the
${\cal S}$-pair $(W, E)$ if it is regular, equidimensional, and has
normal crossings with the divisor $ \cup_{i=1}^{m} H_{i} $. We
define the {\it transformation} of
$(W,E)$ with (permissible) center $C$, denoted by
$(W,E)\longleftarrow (W_{1},E_{1}) $
  to be the pair where $W_1$ is the blowing-up of $W$ with center $C$ and
$ E_{1}=\{H'_{1},\ldots,H'_{m},H_{m+1}\}$, with $ H'_{i}$
the strict
transform of $ H_{i} $, and $ H_{m+1}$ the exceptional
hypersurface in $ W_{1} $. This is again an ${\cal S}$-pair.
\smallskip

(d) An {\it idealistic triple} (or {\it id-triple}) in ${\cal S}$
is a system
$$(1.1.1) \quad {\cal T}=(W, {\cal I}, E)$$
where $(W,E) $ is an ${\cal S}$-pair (as above)  and ${\cal I}$ is
a coherent sheaf of ${\cal O}_W$-ideals
  (a $W$-ideal, or just an ideal, if $W$ is clear, for short) such
that the stalk ${\cal I}_w$ is not zero, for all $w \in W$.

Given such an id-triple ${\cal T}$, let $C$ be a permissible center for
$(W, E)$. If also
$C \subseteq V({\cal I})$, we say that $C$ is {\it permissible} for
${\cal T}$ and we define the {\it transform} of  ${\cal T}$
relative to
$C $ to be the triple
${\cal T}_1:= (W_1,{\cal I}_1, E_1)$, where
$(W_1, E_1)$ is the transform of $(W, E)$ (with center $C$) and
${\cal I}_1={\cal I}{\cal O}_{W_1}$. This is again an id-triple in
${\cal S}$.
\smallskip

Similarly, if $f: W' \to W$ is a smooth morphism, of if $f$ is
obtained by base change from $k$ to a field extension $K$ (recall
that $W$ is a scheme over a field $k$), or if $f$ is a localization
morphism, then an id-triple ${\cal T}$ as above induces an id-triple

${\cal T}':= (W', E', {\cal I}')$, where
$E'= \{ f^{-1}(H_i) : i =1, \ldots ,\, m \}$ and ${\cal I}'={\cal
I}{\cal O}_{W'}$.
\smallskip

(e) Let $X$ be a (possibly empty) closed subscheme of a regular
scheme $W$,
$E=\{H_1, \ldots, H_m \}$ is a collection of regular hypersurfaces
of $W$,
$D=\bigcup _{i=1} ^m H_i$. We say that $X$ has normal crossings with
$E$ if
near each point $w$ of $X_s \cup D$, this
subscheme of $W$ is defined by an ideal
of the form
$( x_1,\ldots, x_m ).
  x_{i_1}.\ldots. x_{i_j} {\cal O}_{W_s,w}  $
where $x_1, \ldots, x_n$   is a regular system of parameters of
${\cal O}_{W_s,w}$, $m \geq n$. If $X$ is empty, we say that $E$ has
normal crossings.
\smallskip

If $ S $ and  $ H $ are regular
  closed subschemes of $ W $ (their defining $W$-ideals being $I(S)$
and $I(H)$ respectively), we say that $S$ and $H$
have a {\it transversal
intersection } at a common point $ x $ if there is a
regular system of parameters $x_1,\ldots, x_n$ of ${\cal O}_{W,x}$
such that:
(1)  $ I(H)_x \subset {\cal O}_{W,x} $ is the ideal generated by $
x_1,...,x_r $,
(2) $ I(S)_x \subset {\cal O}_{W,x}$  is the ideal generated by $
x_s,...,x_n $,
(3) $ r < s $.
\medskip

{\bf (1.2) Definition.} A {\it strong principalization algorithm}
on an allowable collection ${\cal S}$ is an assignment which
associates to each id-triple (in ${\cal S}$)
${\cal T}_0:= (W_0, {\cal I}_0, E_0)$ a
sequence of functions
$h_0, \ldots ,h_r $ ($r$ depends on ${\cal T}_0 $) , which satisfies the
following requirements:
\smallskip

(i) These functions take values in a certain totally ordered set
$I^{(d)}$, $d=\dim W$ (this set depends on $d$ only),
each $h_i$ is
upper-semicontinuous and takes finitely many values.
\smallskip

(ii) The domain of the function $h_0$ is $W_0$ , let $\Max h_0$ be
the maximum value of
$h_0$, and
$C_0=\{ x \in X_0: h_0(x) = \Max h_0 \}$ (we shall also write
$ C_0=\MaxBar h_0 $). Then $C_0$ must be a a
permissible center for
${\cal T}_0$. If $r > 0$, take the transform
${\cal T}_1:= (W_1, E_1, {\cal I}_1)$ of ${\cal T}_0:= (W_0, {\cal
I}_0, E_0)$ with center $ C_0= \MaxBar h_0 $. Then, the
domain of
$h_1$ is $W_1$
and, if
$C_1:= \MaxBar h_1 $, this must be a permissible
center. If $r >1$, take the transform of
${\cal T}_1:= (W_1, E_1, {\cal I}_1)$ with center $C_1$, and so on.
Finally it is required that
$h_r:W_r \to I^{(n)}$ be a constant function and, in
${\cal T}_r:= (W_r, {\cal I}_r, E_r)$, the sheaf of ideals
${\cal I}_r$ be locally principal, with support the union of the
hypersurfaces in $E_r$. (Note that, for each $i$, the hypersurfaces
in $E_i$ have normal crossings, also (by the assumed semicontinuity)
  for each element $ \alpha \in I^{(n)}$, the set
$ F_{i \alpha}=\{\xi \in W_i\,{\rm such \,that} \,
h_i(\xi) \geq \alpha\} $ is closed).

\smallskip

Let
$\Pi_i:(W_i, {\cal I}_i, E_i, ) \to (W_{i-1}, {\cal I}_{i-1}, E_{i-1})$
be the transformation with center $C_{i-1}$, $i=1, \ldots , r$; the
sequence
$$ (1.2.1) \qquad (h_i, \Pi_i:(W_i, {\cal I}_i,  E_i) \to (W_{i-1},
{\cal I}_{i-1}, E_{i-1}), C_{i-1})), \quad i=1,\ldots , r$$
described in (ii) will be called the {\it principalization
sequence}, or {\it p-sequence} corresponding to the id-triple
${\cal T}_0=(W_0, {\cal I}, E_0) $. The function $h_i$ is called
the $i$-th resolution function of the id-triple ${\cal T}_0$.
\smallskip

Furthermore, we require the following conditions:
\smallskip

(iii) If $ X_0=V({\cal I}_0) $, $E_0=\emptyset$  and $\Pi:W_{r} \to
W_{0} $ is the composition
$\Pi=\Pi_r \ldots \Pi_1$,  then  $\Pi:W_{r} \to W_0 $ induces an
isomorphism
$W_{r}-E_{r} \cong W_{0}-X_0 $.
\smallskip

(iv) If $ \xi\in W_i, i=0,\ldots,r-1 $, and if $ \xi \not\in C_{i}$ then
$ h_{i}(\xi)=h_{i+1}(\xi') $
where $\xi '$ is the only point of $W_{i+1}$ corresponding to $\xi$.
\smallskip

(v) $\Max {h_{0}}>\Max{h_{1}}> \cdots >\Max{h_{r-1}} $
\smallskip

(vi) The dimension of each center $C_i$ is determined by the  value
$ \Max {h_{i}}$.
\smallskip

(vii) Let, in (1.2.1), $W'\subseteq W_0$
be an open set. Set $W'_i = (\Pi_i \ldots\Pi_1)^{-1}(W')$,
$C'_i=C_i \cap W'_i$, $E'_i=\{H \cap W'_i: H \in E_r \}$, finally
$h'_i={h_i}_{|W'_i}$. Thus, for all $i$,
there is a naturally defined triple
${\cal T}'_i=(h'_i, \Pi' _i , C'_i)$ (where
${\Pi}'_i:{\cal T}'_i \to {\cal T}'_{i-1}$ is induced by $\Pi _i$).
Then, after neglecting those $\Pi'_i$ inducing
isomorphisms,
the resulting sequence is the p-sequence of the
id-triple ${\cal T}'_0$.
\smallskip

(viii) If ${\cal I}_0$ defines a regular , pure dimensional
subscheme $X_0$ of $W_0$ and $E_0 = \emptyset$, then the
function $ h_{0}$, restricted to $X_0$, is constant.
\medskip

{\bf (1.3) Remarks.} (a) Probably, the condition that looks more
technical and unnatural is (viii). Let us briefly explain some
reasons to include it. First, intuitively the function $h_0$
measures the ``complexity'' of the stalks
${\cal I}_w$, $w \in W_0$. In other words, $h_0(w)=h_0(w')$ should
indicate that the ideals ${\cal I}_w$ and ${\cal I}_{w'}$ are
``equally complicated''. Now, if
$X_0:=V({\cal I})$ is a regular subscheme of $W_0$, then at each
point $w \in W$ the ideal ${\cal I}_w \subset {\cal O}_{W_0,w}$ is
defined by a subset of a regular system of parameters of ${\cal
O}_{W_0,w}$. So, it is reasonable to consider all these ideals as
being ``equally complex'', i.e. $h_0$ is constant on $X_0$. But
there is another important reason for the requirement (viii). It is
known that there is a close connection between the problem of
principalization of ideals and that of embedded desingularization
of varieties. In fact, in (1.6) we shall see how a strong
principalization algorithm (easily) implies one for
desingularization of embedded varieties. Property (viii) will play
an essential role in the proof of this result.
\smallskip

\smallskip

(b) Consider an id-triple (in ${\cal S}$) ${\cal T}=
(W,{\cal I},E)$ (as in 1.1.1) and a  morphism $\phi:W' \to W$ which
is of one of the following types: (i) it comes from an arbitrary
change of the base field, or (ii) it is smooth, or (iii) it is a
localization morphism. Let $W'$ be the induced id-triple ${\cal T}'$
on $W'$ (see  1.1(e)). Then it is easily seen, by using basic
properties of the blowing-up process, that the p-sequence of ${\cal
T}$ (1.2.1) induces a sequence of permissible transformations,
starting at ${\cal T}'$, which is a principalization of ${\cal T}'$.
In fact, the center $C'_0= \MaxBar h'_0$, where $h'_0$
is induced by $h_0$ via the morphism $\phi$, coincides with
${\phi}^{-1}(C_0)$, and the blowing-up of $W'$ with center
$C'_0$ can be identified with the fiber product of $W_1$ and $W'$
over $W_0$, and so on.
\smallskip
Now we shall introduce a notion of algorithmic embedded desingularization
of varieties. First, we need some preliminary definitions.
\medskip

{\bf (1.4)} A {\it couple} in ${\cal S}$ is an ordered pair ${\cal
P} = (X,W)$
where $W$ is a regular, equidimensional scheme which is a member of
${\cal S}$ and $ X $ is a closed subscheme of $ W $.
\smallskip

If ${\cal P} = (X, W)$ is a couple in ${\cal S}$ and $C \subseteq W$ is a
closed and regular subscheme and we take the blowing-up $\Pi: W_1
\to W$ of
$W$ with center  $C$, then $W_1$ is again regular, equidimensional
and, by (2),  a member of  ${\cal S}$. We may
take the strict transform $X_1$ of $X$
in $W$ i.e., $X_1$ is the scheme-theoretic closure of  $\Pi
^{-1}(X-C)$ in
$W$. If $C \subseteq X$, then $X_1$ can be identified to the
blowing-up of
$X$ with center $C$.
The
ordered pair ${\cal P}_1=(X_1, W_1)$ is a new couple in ${\cal S}$,
called the {\it transform} of ${\cal P}$.
\medskip

{\bf (1.5) Definition.} A {\it desingularization (or resolution)
algorithm} on an allowable collection ${\cal S}$ is an assignment
which
associates to each couple ${\cal P}=(X_0, W_0)$ in ${\cal S}$,
where $X_0$ is reduced and equidimensional, a
sequence of
functions
$f_0, \ldots ,f_s $ ($s$ depends on ${\cal P}$) , which satisfies the
following requirements:
\smallskip

(i) These functions take values in a certain totally ordered set
${\Lambda}^{(n)}$, $n=\dim W$ (this set depends on $n$ only),
each $f_i$ is
upper-semicontinuous and takes on finitely many values.
\smallskip

(ii) The function $f_0$ has $W_0$ as domain, let $\Max f_0$ be
the maximum value of
$f_0$, and
$C_0=\{ x \in X_0: f(x) = \Max f_0 \}$ (we shall also denote
$ C_0= \MaxBar f_0 $). Then $C_0$ must be a regular
(necessarily closed) equidimensional subscheme of $W$.
\smallskip

(iii) If $f_0, \ldots, f_i$ have been defined, with
$C_j=\{ x \in W_j: f(x) = \Max f_j \}$ , $j=0, \ldots, i$ being a
regular,
equidimensional subscheme of $W_j$; let $(X_{i+1},W_{i+1}) \to
(X_{i},W_{i})$
be the transform of
$(X_{i},W_{i}) $ with center $ C_i $. If $\Pi=\Pi_1 \ldots \Pi_i$ set
$ E_{i+1}=\{H_1,H_2,\cdots,H_{i+1} \}$
the set of hypersurfaces in $W_{i+1}$ ,where each $H_j$ is the strict
transform
of the exceptional locus of $ \Pi_j$;then $E_{i+1}$ is a set
of smooth hypersurfaces having only normal crossings.
Then the domain of $f_{i+1}$ is $W_{i+1}$, and
$C_{i+1}$ (the set
of points where it reaches its maximum) must be a regular,
equidimensional
subscheme of $W_{i+1}$ having only normal crossings with  $E_{i+1}$.
Since $ C_i $ is the closed set where $ f_i $ takes maximum value
($\Max f_i $), we
shall sometimes write $ C_i:= \MaxBar f_i $.
\smallskip

(iv) The strict transform $X_s$ of $ X_0 $ to $ W_{s} $ is
regular and has normal crossings with $E_s$. This means that,
letting $D$ denote the union of all the hypersurfaces in $E_s$ (with
reduced structure), then at each point $w$ of $X_s \cup D$, this
subscheme of $W_s$ is locally defined by an ideal
of the form
$( x_1,\ldots, x_m ).
  x_{i_1}.\ldots. x_{i_j} {\cal O}_{W_s,w}  $
where $x_1, \ldots, x_n$   is a regular system of parameters of
${\cal O}_{W_s,w}$.
\smallskip

The sequence
$$ (1.5.1) \qquad (f_i, \Pi_i:(X_i,W_i) \to (X_{i-1},W_{i-1},
C_{i-1})), \quad
i=1,\ldots , s$$
described above will be called the {\it resolution sequence}
corresponding to the couple ${\cal P}=(X,W) := (X_0, W_0)$.
\smallskip

As announced before, next we explain how the presence of a strong
principalization algorithm (Definition 1.2) implies the existence of
one of desingularization (Definition 1.5).
\medskip

{\bf 1.6 (From principalization to desingularization) } Assume a
strong principalization algorithm on an allowable collection ${\cal
S}$ (1.2) is given. We obtain a desingularization algorithm on $\cal
S$ as follows.
Consider an ${\cal S}$-couple $(X_0,W_0)$, where
$ X_0 $ is a reduced and pure dimensional subscheme of the
$d$-dimensional regular scheme $ W_0 $,
defined by a coherent sheaf of ideals $ I(X_0) \subset {\cal
O}_{W_0} $. Set $ {\cal I}_0=I(X_0) $ and  take the id-triple $
{\cal T}_0=(W_{0},{\cal I}_{0},\emptyset) $ and the p-sequence
(1.2.1) corresponding to $  {\cal T}_0 $. Recall that,
by property 1.2(viii), the function $ h_{0}$
is constant , say equal to $\lambda \in I^{(d)}$, along the
non-empty open set
${\rm Reg}(X_{0})$ of regular points of $X_0$.
     By 1.2(iv) and (v) we see that there must be
     a unique index $ s \leq r-1 $  such that
     $\Max{h_{s}} = \lambda $.
Let $X_i$ be the strict transform of $X_{i-1}$ via
$\Pi _i : W_i \to W_{i-1}$ (the blowing-up with center $C_{i-1}$),
for $i \geq 1$.
Applying once again 1.2(iv) we see
that the naturally induced morphism $p: X_{s} \to X_0$ restricts to
an isomorphism $p^{-1}({\rm Reg}(X_{0})) \to {\rm Reg}(X_{0})$.
Moreover,
$X_{s}$ must be included in the closed
     set $ C_{s}=\MaxBar{h_{s}} $.
But (by (ii) of 1.2)  the centers provided by the algorithm of
principalization
    have normal crossings with hypersurfaces in $ E_{s} $.
     Finally, the smoothness and pure dimensionality of $C_{s}$
and 1.2(vi) imply that $X_{s}$ must be a union of
     connected components of $C_{s}$. Hence $ X_{s} $ is
     smooth and has only normal crossing with $ E_{s}$. This shows
that the first $ s $ steps define our
embedded desingularization. From this,
$f_i=h_i$, $1=0, \ldots, s$ satisfies the requirements of 1.5 and
we have an algorithm of desingularization for the couple
$(X_0,W_0)$.
\smallskip

The index $s$ we obtained is called the {\it resolution index} of
the couple
$(X_0,W_0)$.
\medskip

{\bf (1.7) Remark.} The result of 1.6 is not new, it appears in
[12], although there one works with a specific principalization
algorithm. Originally, the goal of the present  paper was to study
families of embedded schemes and criteria that insure that all the
members can be simultaneously desingularized, using an algorithmic
resolution process fixed beforehand. The discussion of 1.6 suggests
that it is more more natural to emphasize the study of families of
ideals and of ways to simultaneously principalize its members (by
means of a given algorithm.) It seems that to succeed one must impose
on the algorithm some further conditions (aside from those listed
in 1.2). So, next (after a few explanatory remarks) we shall study
some basic notions about families
of ideals and discuss the ``extra conditions''.
\medskip

{\bf (1.8) Localizations.} Given a morphism $f:W \to T$ of schemes
and a point $t \in T$, consider the canonical morphism
$T':={\rm Spec}\,{\cal O}_{T,t} \to T$ sending the closed point  
$t'$ of $T'$ to
$t$. The resulting morphism
$f\{t\}:W' \to T'$ obtained by base change is called the {\it
localization} of $f$ at $t$. Often we shall write $T':=T\{t\}$ and
$W':=W\{t\}$ and, if $T$ and $f$ are clear from the context,
$W\{t\}$
itself will be called the localization of $W$ at $t$. Note that
there is a canonical isomorphism of residue fields $k(t)=k(t')$ and
the fiber of $f$ at $t$ can be identified to that of $f'$ at $t'$.
Moreover, the natural morphism $j_t:W\{t\} \to W$ is a localization
morphism (cf. 1.1(a).)
\medskip

{\bf (1.9) Remark on the notation.} (a) Let $\pi:W \to T$ be a
smooth morphism of regular noetherian schemes, $t$ a point of $T$.
As usual, the fiber $\pi$ at $t$ is the $k(t)$-scheme $W^{(t)} \to
{\rm Spec}\,(k(t))$ obtained from $\pi$ by the base change
${\rm Spec}\,(k(t)) \to T $ (the canonical morphism whose image is
$t$). Similarly, if instead we use the canonical morphism
${\rm Spec}\,({\overline {k(t)}}) \to T $, where ${\overline
{k(t)}}$ denotes an algebraic closure of $k(t)$, we obtain the
geometric fiber at $t$ (this is, up to isomorphism, independent of
the choice of the algebraic closure). Sometimes we denote this
geometric fiber by $W^{(\bar t)}$.

(b) We have a canonical morphism $j_t:W^{(t)} \to W$, if $t$ is a
closed point of $T$ then $j_t$ is a closed embedding. In general, if
$Z$ is a closed subscheme of $W$ the symbol $Z \cap W^{(t)}$ (or $Z
\cap {\pi}^{-1}(t)$) will denote
${j_t}^{-1}(Z)$ (a closed subscheme of $W^{(t)}$).
If $t$ is closed, then via the
closed embedding $j_t$ this becomes the usual intersection of
subschemes of
$W$. Note that always the fiber $W^{(t)}$ can be identified with
the closed
subscheme ${\pi \{t\}}^{-1}(t')$ of the localization $W\{t\}$ (cf.
1.8). Then, if $Z\{t\}:= j_t^{-1} (Z) \subset  W\{t\}$ (notation of
1.8),
$Z \cap W^{(t)}:=j_t^{-1}(Z)$ gets identified to the actual
scheme-theoretic
intersection $Z':= Z\{t\} \cap {\pi \{t\}}^{-1}(t')$ (inside $W
\{t\}$). Sometimes, when the meaning is clear from the context, $Z'$
will be also denoted by $Z \cap W^{(t)}$.
\medskip

{\bf (1.10) Definitions.} (a) A {\it family of ideals} is a triple:
$$(1.10.1)\quad {\cal G}=(\pi:W \to T, {\cal I}, E)$$
where $\pi$ is a smooth morphism and $(W, {\cal I}, E)$ is an
id-triple (1.1.1), satisfying the following condition.  If $t$ is a
point (or a geometric point) of $T$, let
$W^{(t)} = {\pi}^{-1}(t)$, ${\cal I}^{(t)}={\cal I}{\cal
O}_{W^{(t)}}$ and
$E^{(t)}=\{H^{(t)}_{1}, \ldots,H^{(t)}_{m} \}$, where
$E=\{H_{1}, \ldots,H_{m} \}$ and
$H^{(t)}_{i}=H_i \cap W^{(t)}$, for all i (cf. 1.9 for the meaning
of the intersection symbol). It is required that, for
all $t \in T$,
${\cal G}(t):=(W^{(t)}, {\cal I}^{(t)}, E^{(t)})$ be an id-triple,
i.e., the $H^{(t)}_{i}$'s are hypersurfaces of $W^{(t)}$ with normal
crossings only.

The id-triple ${\cal G}(t)$
will be called {\it the fiber at $t$}, we'll also say that ${\cal
G}(t)$ is a member of the family ${\cal G}$.

If $(W,E)$ and all the fibers are  ${\cal S}$-pairs (where ${\cal
S}$ is an allowable
collection), we
say that we have an ${\cal S}$-family.

If $t \in T$, consider the {\it localization}
  $W\{t\} \to T\{t\}$ of $\pi$
  at $t$ (see 1.8). As above for the fiber ${\cal G}(t)$, ${\cal G}$
induces a family of ideals
${\cal G}\{t\}=({\pi}\{t\}:W\{t\} \to T\{t\},{\cal
I}\{t\},E\{t\})$, called the {\it localization of ${\cal G}$ at
$t$.} Note that $T\{t\}$ has only one closed point and the fiber
$W^{(t)}$ of $\pi$ at $t$ can be identified to the closed fiber of $
{\pi}\{t\} $.
\smallskip

By the p-sequence of the family 1.10.1 we mean the p-sequence of the
id-triple
$(W,{\cal I},E)$.
\smallskip

(b)  A {\it family of embedded schemes} is
a pair
$$ (1.10.2) \qquad {\cal F} = (j: X \to W, \pi:W \to T) $$
where  $W$ and $T$ are noetherian schemes over a field $k$ of
characteristic zero, with $W$  regular, of pure dimension $n$, $j$
is a closed immersion, $\pi$ is smooth and $p:= \pi j:X
\to T$ is flat. We
write $X^{(t)}:=p^{-1}(t), W^{(t)}={\pi}^{-1}(t)$, for $t \in T$.

If ${\cal S}$
is an allowable collection (see (1.1)(a)), we say that a family of
embedded schemes is an
${\cal S}$-family if $(X,W)$ and all the fibers
$(X_t,W_t)$ are couples in ${\cal S}$.
\smallskip

(c) Given a family of embedded schemes (1.10.2) we have an {\it
associated family} of ideals, given by
$$(1.10.3) \quad (\pi, I(X), \emptyset)$$
where $I(X)$ is the $W$-ideal defining the subscheme $X \subset W$.
\smallskip

(d) Note that if ${\cal G}$ (1.10.1) is a an ${\cal S}$-family of
ideals, then all the fibers are id-triples in $\cal S$.
\medskip

{\bf (1.11) The notion of compatibility.} Given an algorithm of
strong principalization on an allowable collection ${\cal S}$ and
and a family of ideals ${\cal G}$ (1.10.1), we may consider the
id-triple
$({\cal I}, W, E)$ and its associated p-sequence (1.2.1) (with
$({\cal I}, W, E)=({\cal I}_0, W_0, E_0) $. On the other hand, if
$t \in T$, we may consider the triple
${\cal G}(t)=({\cal I}^{(t)}, W^{(t)}, E^{(t)})$ (the fiber at $t$)
and its associated p-sequence, say
$$(1.11.1) \quad
(h^{(t)}_i, \Pi ^{(t)}_i:
(W^{(t)}_i,  {\cal I}^{(t)}_{i},   E^{(t)}_i)
\to
(W^{(t)}_{i-1}, {\cal I}^{(t)}_{i-1}, E^{(t)}_{i-1}),
C^{(t)}_{i-1})), \quad i=1,\ldots , r_t.$$
We would like to compare the values of the functions $h_i$ and
$h^{(t)}_i$. But, in general, beyond $i=0$ there is no way to
identify the domain of
$h^{(t)}_i$ (namely, $W^{(t)}_i$) with  a subset of $W_i$ (the
domain of $h_i$), hence the comparison makes no sense. We want to
introduce a condition that will make this comparison possible.
\smallskip

Given a triple
${\cal T}=(W, {\cal I}, E)$ (in ${\cal S}$) and a closed, regular
subscheme
$\overline{W}$
of $W$, next we define the concept `` $\overline{W}$ is
j-compatible with ${\cal T}$ '' (relative to a given strong
principalization algorithm, when the choice is clear from the context we
won't mention the algorithm).

If $j=0$, this just means that for all $H \in E$, $H$ and
$\overline{W}$ intersect transversally (1.1.e) along a hypersurface
$\overline{H}$ of $\overline{W}$, and the collection
$\overline{E}$
of these $\overline{H}$'s have normal crossings. Note that from this, if
$  \overline{{\cal I}}={\cal I}{\cal O}_{W'}$, then
  $\overline{\cal G}= (\overline{W}, \overline{{\cal I}},
\overline{E})$ is a new id-triple in ${\cal S}$. It has a
p-sequence, its elements will be denoted as in (1.2.1), but with
bars on top.

If $j >0$, we require, moreover, that
$C_0 : = \MaxBar h_0$ be transversal to
$\overline{W}$, and
$C_0 \cap \overline{W}=\overline{C_0}$ (where $C_0$ and
$\overline{C_0}$ are the 0-centers in the p-sequences of ${\cal G}$ and
$\overline{\cal G}$ respectively). It is known that this implies that
$\overline{W}_1$, the blowing-up  $\overline{W_1}$  of
$\overline{W}$ with center
$\overline{C_0}$, can be identified with the strict transform of
$\overline{W}$ via $\Pi _1:W_1 \to W$ (the blowing-up with center
$C_0$) and the exceptional locus $\overline{H}$ with $H \cap
\overline{W_1}$, where $H$ is the exceptional locus of $W_1 \to W$.
If $ j >1$, we further require (using this identification):
$C_1$ is transversal to $\overline{W_1}$ and $C_1 \cap W_1 =
\overline{C_1}$. Then,
${\overline{W}}_2$ may be identified to the strict transform of
$W_1$ (contained, up to canonical identification, in $W_1$) to $W_2$
(the blowing-up of $W_1$ with center $C_1$. And so on. If this
process can be repeated $j-1$ times, we say that $\overline{W}$ is
{\it $j$-compatible with ${\cal T}$}.

Thus, in this case,
${\overline{W}}_i$ can be identified to a subscheme of  $W_i$, for
$0 \leq i \leq j$. In particular, it makes sense to ask: is $C_j$
transversal to $\overline{W}_j$?

The most important application of this notion is to the case where we
have a family of ideals ${\cal G}$ (1.10.1) and
$\overline{W} \subset W$ is a fiber of $\pi$,
$\overline{W}=W^{(t)}$, $t \in T$.

Precisely, we say that
${\cal G}$ is
$j$-compatible with the algorithm at $t $ if
${\cal T}\{t\}:=( W\{t\}, {\cal I}\{t\}, E\{t\})$ (the id-triple
coming from the localization at of ${\cal G}$ at $t$) is
$j$-compatible with the closed fiber of the localization
${\pi}\{t\}$ of the morphism $\pi$. Note that if $t$ is a closed
point of $T$, this is equivalent to saying that  $(W,{\cal I},E)$ be
$j$-compatible with the fiber of  $\pi$ at $t$ (a closed subscheme
of $W$). Finally, we say ${\cal G}$ is
$j$-compatible with the algorithm if there is $j$-compatibility at
each $t \in T$. Of course, this forces the index
$r_t$ of (1.9.1) to be greater than or equal to $j$, for any point
$t$ of $T$.
\medskip

{\bf (1.12) Remark.} We make some comments about certain concepts
used in 1.11. Here, $ \pi:W \to T $ is a  morphism of
finite type of regular noetherian schemes, $S \subset W$ is a
regular subscheme.
\smallskip

(a) If a closed point
$ x \in W $ maps to $ t \in T $, then the morphism is smooth
  at $ x $
  if and only if a regular system of parameters
$x_1,\ldots, x_r$ of ${\cal O}_{T,t}$ can be extended to
a regular system of parameters ${\cal
O}_{W,x}$.

If $\pi$ is smooth and $S$ (containing $x$) is as above,  then the
induced morphism  ${\eta}: S \to T $
is also smooth  at $ x$ iff $ S $ and $ H:={\pi}^{-1}(t) $ intersect
transversally at $ x $.
\smallskip

(b) If $ \pi$ and the induced morphism  $ {\eta}:S \to T $ are both
smooth,
and  $ W' \to W $ is the monoidal transformation with center $ S $, then
the composite morphism $ \pi':W' \to T $ is also smooth. We sketch
the proof of this known fact.
  The local ring
${\cal O}_{W',x'}$ at a closed point $ x' \in W' $ mapping to $ x \in W $
is a localization of ${\cal O}_{ W,x } [v_s,...,v_n ]$ at a maximal
ideal.
(We may assume that there is an index $j_o$ so that $ x_{j_0}v_k=x_k $
for $ k= s,s+1,...,n $.) One can check now from 1.1(e), (3), that the
regular system of parameters $x_1,\ldots, x_r $ of $ {\cal O}_{T,t} $ can
be extended to a regular system of parameters of $ {\cal O}_{ W',x'} $.
Then use 1.1(e).
\smallskip

{\bf (1.13) Good principalization algorithms.} A strong principalization
algorithm is {\it good} if
it has the
following properties.
\smallskip

(1) {\it Compatibility with change of the base field.}
Assume that, in the id-triple ${\cal T}_0=(W_0, {\cal I}_0,  E_0)$
(in ${\cal S}$), where $W_0$  is a scheme over a field $k$. Let $ k
\subset K $ be a field
extension. Consider the resolution sequence (1.2.1) of ${\cal T}_0$.
Then (see 1.3(b)) there is an induced
pricipalization sequence over $ K $, say
${\cal T}'_0 \leftarrow \ldots \leftarrow {\cal T}'_r$. This must
coincide with the p-sequence that the algorithm attaches to ${\cal
T}'_0$.
\smallskip

(2) Consider the p-sequence 1.2.1 of an i-triple ${\cal T}_0$ and  ,
for  $\alpha \in I^{(n)}$,
$S_{i,{\alpha}}= \{\xi \in
W_i: h_i(\xi) = \alpha\}$. Then $S_{i,{\alpha}}$ (a locally closed
subset of $W_i$, regarded as a reduced subscheme) is regular and
pure dimensional.
\smallskip

(3) {\it Compatibility with $pr_1:W \times {\bf A}^m \to W $}.
There are order-preserving injective functions
$$ (1.13.1) \quad \lambda_{n,m}:I^{(n)} \to I^{(n+m)}$$
defined for integers $n>0$,$m>0$, with the following property.
Let ${\cal T}_0=(W_0, {\cal I}_0,  E_0)$ be any id-triple in ${\cal
S}$ (where  dim $W_0=n$) having
p-sequence (1.2.1). Let $W'_0:=W_0 \times {\bf A}^m$ and
     $ p_0:= W'_0 \to W_0$
be the first projection. Consider the corresponding induced
sequence, obtained by pull-back,
$ (W_0, {\cal I}_0, E_0)\leftarrow (W_1, E_1, {\cal I}_1)
\leftarrow \ldots  \leftarrow (W_r, E_r, {\cal I}_r)$ (see 1.1).
(Note that for
all $ i $ there is a morphism $ p_i: W'_i \to W_i $). Then this coincides
with the p-sequence
$(h'_i, {\Pi}'_i, C'_{i-1})$
that the algorithm attaches to
$(W'_0, E', {\cal I}')$; moreover for each $i$ and
$w'_i \in W'_i$, $h'_i(w'_i)=\lambda_{n,m}h_i p_i (w'_i)$
\smallskip

(4) {\it Compatibility with localizations.} Let, in (1.2.1), $\phi:
W' \to W$ be a localization morphism
(1.1.a). By  1.3 (b) , the p-sequence of ${\cal T}=(W_0,{\cal I}_0,
E_0)$ induces a principalization sequence of the id-triple ${\cal
T}'$, obtained from ${\cal T}$ by pull-back via $\phi$. Then
neglecting those steps where (using the notation of 1.2.1)
$C_i \cap W\{t\}=\emptyset$,
the resulting sequence is the p-sequence of the
id-triple ${\cal T}'_0$.
\smallskip

(5) Consider an ${\cal S}$-family (1.10.1) ${\cal G}$ which is
$j$-compatible
with the algorithm, with p-sequence (1.2.1). Let
(1.11.1) be the sequence for the fiber at $t \in T$.
Then, we require (using
the notation of (2)): if $x \in S_{j,{\alpha}} \cap
(W_j)_t$, then
${h_j}(x) \leq \lambda_{n,m}(h_j ^{(t)}(x))$
(with $ \lambda_{n,m}:I^{(n)} \to I^{(n+m)}$ as in part (3)), and
equality holds if and only if
the intersection is transversal at $x$ (see 1.1(e) and the next
Remark 1.14).
\medskip

{\bf (1.14) Remark.}
In writing the inequality in 1.13(5), we are making some
identifications, as follows. Consider the localization ${\cal
G}\{t\}$ of the family  at $t$. By (4), the p-sequence of ${\cal G}$
induces that of ${\cal G}\{t\}$, in particular the j-th function
$h_j$  of ${\cal G}$ induces a similar function $h_j'$ for ${\cal
G}\{t\}$, with domain $W_j\{t\}$. The fiber $W_j^{(t)} \subset W_j$
corresponds to the closed fiber $F$ of the natural morphism
$W_j\{t\} \to T_j\{t\}$, and (by $j$-compatibility), $W_j^{(t)}$
can be identified to $F$. The inequality above really takes place on
$W_j\{t\}$, involving ${h_j}'$ and the function (defined on $F$)
naturally corresponding to $h_j^{(t)}$. Now, $S_{j,\alpha}$
corresponds to a closed subscheme  $S'_{j,\alpha}$ of $W_j^{(t)}$,
the transversality of the intersection really refers to $F \cap
S'_{j,\alpha}$ (in $W_j^{(t)}$). If $t$ is a closed point, then
$W_j^{(t)}$ is a closed subscheme of $W_j$ and there is a simpler
interpretation, working inside $W_j$ and not using explicitly
localization.
\smallskip

In 1.13 (5) (a property which will play a central role throughout
this paper)
in general we shall simply write ${h_j}(x) \leq h_j ^{(t)}(x)$
identifying the
right-hand-side and its image in $ I^{(n)}$ via $\lambda_{n,m}$.
Note that
the transversality at $ x $ ( equivalent to the equality) implies
that, if $ x \in W_t$ and $ S_{j,\alpha}$
is the stratum containing
$ x $, then $x$ must be a regular point of both $ W^{(t)} \cap
S_{j,\alpha}$ and
$ S_{j,\alpha}$.
\medskip

{\bf (1.15)} Now we present an example of a good
principalization algorithm.  Let ${\cal S}_0$ be the class of
regular, equidimensional schemes $ W $
containing a field, say $k$, of characteristic zero (which may
vary), satisfying:

(i) If $ W $ is an n-dimensional $k$-scheme in ${\cal S}_0$ , then
it has a
finite
affine open covering $\{{\cal U}_j\}$, $j \in J$ where, for each
$j$, ${\cal U}_j \approx
Spec(R_j)$, with $R_j$ being a noetherian, regular $k$-algebra, and such
that $Der _k(R_j)$ is a finite projective $R_j$-module, locally of
rank n.

(ii) if $ M $ is a maximal ideal in $ R_j $ then $ R_j/M $ is
algebraic over $ k $.

Under condition ii) $ k $ is a quasi-coefficient field at the
localization
at any closed point in the sense of [21], page 274.

Then we have:

\proclaim
(1.16) Theorem. Let ${\cal S}_0$ be the class of 1.15.
Then ${\cal S}_0$ is an allowable collection and there  exists a good
principalization algorithm on ${\cal S}_{0}$.

{\bf (1.17)} It is easy to verify  that ${\cal S}_{0}$ is an
allowable collection. We won't give a complete proof proof of the
theorem. The bulk of the necessary work can be seen in [11] and
[10]. In these references, there is detailed account of the
construction of an algorithm to resolve {\it basic objects} (see
5.1). From this, the existence of a principalization algorithm
easily
follows (see [12], or 5.5 of the present paper) All the required
properties of a good algorithm are covered in [11] and [10],
except 1.13(4) and 1.13(5). These are proved in 1.18 and chapter~6
respectively.
Moreover, in chapter~5 we shall review some material from [11],
necessary for these proofs., and we indicate how an algorithm to
resolve basic objects implies a principalization algorithm for
triples.

Actually, in [11] and [10] one works with schemes of finite
type over  fields. However, the constructions and proofs of these
papers are valid over ${\cal S}_0$. In fact,   all what we need is
the following constructions,
available in ${\cal S}_0$, and the validity of Property D, stated
below.

If we fix a ring $ R_j $ as in i) of 1.15, say $ R $, we  introduce an
operator on
the set of all ideals of $R$ by setting  $ \Delta=Der _k (R)$ and
defining, for an ideal $J$ of $R$,  $\Delta J$ as the ideal of $R$
generated  $J$
and the set $ \{\delta(f) / f \in J,\delta \in \Delta \} $.
We refer here to [21], appendix 40,
particularly theorems 99 and 102. It can be checked, using the equivalent
conditions (3) and (4) in Th. 99, that the class is closed by monoidal
transformations.

Note that  any smooth scheme over a field of characteristic zero,
as well as the spectrum of the completion or henselization of a
local ring thereof at a closed  closed point is in ${\cal S}_0$.

Now, the following property holds:
\smallskip

PROPERTY D: If $W=Spec (A) \in {\cal S}_0$, $P \in W$ and $J$ is any
ideal of $A$, then the order of $ J A_P $ in the local regular ring
$A_P$ is $ \geq b$
  iff $ P $ contains the ideal $ \Delta^{b-1}(J) $.
\smallskip

To check this property, let $A$ be as above, $ R $ the localization
of $A$
at a maximal ideal of $A$ and
  $R'$  the completion of $R$. Then, the residue field of $R'$ is
a finite extension of $k'$ of $k$
and $ Der_{k}(R) $ induces $
Der_{k'}(R') $
over $ R' $ (see Th. 99,(4)). By Theorem 102 of [21] we know
that $R$
is excellent, so $ R \to R'$ is faithfully flat with regular
fibers.
Assume first that $P$ is a regular prime ideal in $R $. It induces
a regular prime
ideal in
$ R' $ and it is easy to see that the Property D holds for
regular
primes by checking that it holds at the completion, which is a ring of
formal power series.

Now let $ P
\subset R$  be an arbitrary prime ideal. First note that it
suffices to consider the case in
which $J $
is principal. In this case, set $ \bar{R}=R/J$; the order of $ J $
at $P $
is the multiplicity of the ring $ \bar{R}$ at the prime ideal, say $
\bar{P}$, induced by $ P $.

On the one hand the multiplicity is an upper-semicontinuous function, and
on the other hand any prime ideal $ P \subset R $ is the intersection of
all prime ideals of height $ n-1$. In fact if $ dim R/P= n-h $ and
$ f \in
R $ is not in $ P $, one can find a prime ideal of height $ n-1 $
containing $ P $ but not $ f $ (set $ \bar{f}_1 \in R/P$ as the
class of $
f $, extend to
$ \bar{f}_1 ,\bar{f}_2 ,...\bar{f}_{n-h} \in R/P$ a system of parameters,
and now take $ Q \subset R $ by lifting a minimal prime ideal containing
$ <\bar{f}_2 ,...\bar{f}_{n-h}> \subset R/P$).

  Hence it suffices to assume that $ P $ is a
prime ideal defining a curve. We finally follow a trick of Hironaka
[15]: we reduce to the case where $ Q $ is regular by desingularizing
the curve
by means of a composition of quadratic transformations.
In order to make this reduction possible, we view here $ Der_k (W)$ as a
coherent sheaf over $ W $, and if
$\pi: W' \to W $ is a monoidal transformation, and $y \in W' $ is
isomorphic to $ x \in W $ via $ \pi$, we identify the localizations of $
Der_k (W_1)$ at $y$ with that of
$ Der_k (W)$ at $x$. Finally note the operator $ \Delta $ in
property D is
defined in terms of $ Der_k (W)$. It suffices now to take $\pi: W'
\to W $
as an embedded desingularization of the curve defined by $ P $, and
set $ x
\in W $ and $ y \in W' $ as the localizations at $ P $ and at $ Q $
respectively.
\medskip

{\bf (1.18)} We note now that property 1.13(4), for this algorithm,
follows
from property~D.
Given a localization morphism $W'\longrightarrow W$ and a point in $W'$,
there exists a neighborhood such that the morphism is of type
$Spec(S^{-1}A)\longrightarrow Spec(A)$, where $S$ is some
multiplicatively closed set in $A$.
Now 1.13(4) follows from the fact that property~D is stable by
localization at any multiplicatively closed set.
\bigskip

{\bf 2. EQUISOLVABILITY AND STATEMENT OF THE MAIN THEOREM.}
\medskip

Throughout this section, we work within a given
allowable collection ${\cal S}$ where a strong principalization
algorithm has been defined (1.2).
We want to introduce conditions that insure that all the different
members of
a family of ideals ${\cal G}= (\pi:W \to T, {\cal I}, E)$ (1.10.1)
can be simultaneously principalized, by using the algorithm. When
this happens, we shall say that the family is {\it equisolvable} (to
avoid the use of very complicated words, such as
``equiprincipalizable''.) We propose two definitions (2.2. and 2.3),
which turn out to be equivalent (Theorem 2.4, to be proved in the
next chapter). The first one does not explicitly involve the fibers
of the family, but it requires, essentially, that the centers $C_i$
that appear in the p-sequence (1.2.1) that the algorithm associates
to the id-triple $(W,{\cal I},E)$ be ``evenly spread'' over the
parameter space $T$. In the second condition, a ``numerical''
invariant is associated to the different points $t \in T$ (this
invariant in defined in terms of the p-sequences of the fibers); it
is required that it be constant along $T$. Both approaches have
their advantages, depending on the situation (e.g., see 3.7). In
Chapter 4 we shall see that given an arbitrary family (1.2.1) of
ideals, it is possible to naturally stratify the parameter space
$T$ as  a union of locally closed sets so that, along each one, the
restriction of the family is equisolvable.
\smallskip

Finally recall (1.6) that the given strong principalization
algorithm induces an {\it associated desingularization algorithm}
for couples $(X,W)$ in ${\cal S}$. We shall say that a family of
embedded schemes (1.10.2) is {\it equisolvable} (relative to our
algorithm) if the associated family of ideals (1.10.3) is
equisolvable. If the family (1.10.2) is such that all the fibers
$X_t$ are reduced, then the desingularization sequence (1.5.1) that
the associated algorithm assigns to the couple $(X,W)$ induces on
each fiber the resolution sequence that corresponds to that fiber,
that is it has the property to be expected of a good notion
of simultaneous resolution.

But this definition of equisolvability applies also to the case
where some (or all the) fibers $X_t$ are non-reduced. So, we have a
(we hope, reasonable) notion of equiresolution for families of
embedded schemes where some fibers may be non-reduced. Since in many
geometric problems it is unavoidable the presence of non-reduced
fibers in families of schemes, this seems to be an important
feature.
\medskip

{\bf (2.1)} We introduce the following condition on an ${\cal S}$-
family of ideals ${\cal G}$ (1.10.1). Let (1.2.1) be the p-sequence
  of
$  (W_0, {\cal I}_0, E_0)=(W, {\cal I}, E)$. Note that we get, by
composition,
morphisms $\pi
_{i}: W_i \to T, p_i:X_i \to T$, which induce morphisms $\rho _i: C_i \to
T$, $i=0, \ldots, r$.
\medskip

{\bf Condition (AE).} The morphism $\rho _i: C_i \to T$ is smooth,
proper and
  surjective, for $i=0, \ldots, r$, where $r$ denotes the last index
in the principalization sequence 1.2.1.
\medskip

{\bf (2.2)} Next, given a family ${\cal G}$ (1.10.1), we shall define a
  function $\tau _{{\cal G}}$ from the parameter space
$T$ into a certain totally ordered set ${\Lambda ^{(m)}}$ (which
depends on
$m = $ dimension of the fibers of $\pi: W \to T $ only).
\smallskip

For each $t \in T$ and $i= 0, \ldots, r_t $ ($r_t $ as in 1.11.1), let
$c^{(t)}_i$ be the
number
of connected components of
$C^{(\bar t)}_i := \MaxBar h^{(\bar t)}_i$ (i.e., of the $i$-th  
center of the
geometric fiber at $t$). Then we set:
$$(2.2.1) \quad \tau _{{\cal G}}(t) =
( \Max h^{(t)}_0,c^{(t)}_0, \Max h^{(t)}_1, c^{(t)}_1, \ldots,
\Max h^{(t)}_{r_t},c^{(t)}_{r_t}, \infty, \infty, \ldots).$$

That is, the
values are sequences whose entries are either in {\bf Z}, or in $I^{(m)}$
or (eventually) $\infty$, ordered lexicographically.
\smallskip

{\bf Condition $\tau$.} We require that $\tau _{{\cal G}} (t)$ be
constant, for
all $t
\in T$.
\medskip

Now we state our main theorem.

\proclaim
(2.3) Theorem. Consider a family ${\cal G}=(\pi:W \to T, {\cal I}, E)$
  as in (1.10.1), with $T$
integral and such that, for all i, the morphism $\rho _i : C_i \to
T$ is
proper. Then
it satisfies Condition (AE) if and only if it satisfies Condition $\tau$.
In either case, the principalization sequence (1.2.1) of
$ ( W , {\cal I}, E) $ induces, by taking
fibers,
the principalization sequence of
${\cal G}( t):=(W^{(t)}, {\cal I}^{(t)}, E^{(t)})$
  or of the
geometric
fiber
${\cal G}(\bar t):=(W^{(\bar t)}, {\cal I}^{(\bar t)}, E^{(\bar t)})$,
for all $t \in T$.

The last part of the statement means the following.
The length $r$
of the principalization sequence of $ ( W , {\cal I}, E) $
  agrees with the length $r_t$
of the principalization sequence of any fiber
${\cal G}(\bar t):=(W^{(t)}, {\cal I}^{(t)}, E^{(t)})$, there  is a
natural
identification of $\pi _i ^{-1}(t)$ and $W_i ^{(t)}$ (cf. 1.9)
for any $t \in T$ and, via this identification, the restriction of
$h_i$ to $(\pi _i ^{-1}(t))$ coincides with $h^{(t)}_i$, for $i=0,
\ldots, r$. Similarly with geometric fibers.
\medskip

{\bf (2.4) Definition.} Given a family of ideals as in 2.3, we say
that it is {\it equisolvable} if any of the equivalent conditions of
Theorem 2.3 holds.
\smallskip

Note that the hypothesis in Theorem 2.3 are automatically fulfilled
if $\pi: W \to T$ is proper, in particular projective.
\medskip

{\bf (2.5) Remark.} To prove the Theorem, we shall use the
following results:
\smallskip

(i) Fix a family ${\cal G}=(\pi:W \to T, {\cal I}, E)$ and let
(1.2.1) be the
principalization sequence of $ ( W , {\cal I}, E) $. Fix $j \in
\{0, \ldots, r
\}$, and let $a_{j1}< \ldots < a_{js}$  be the
values of $h_j$ (they are finitely many), then
$S_{jl}=h_j^{-1}(a_{jl})$ is
a regular, locally closed and
equidimensional subscheme of $W_j$. In fact, note first that $h_r$ is
constant along $ W_r $
( where $r$ is the length of the
principalization sequence $(1.2.1)$ of $ ( W , {\cal I}, E) $. Using
now (iv) in (1.2) and condition 2 in
(1.13) we see that, for a given $i \in \{0, \ldots, s \}$, there
must be an index $v$ such that:

a) $\Max h_v=a_{ji}$,

b) If $x \in S_{ji}$, a neighborhood of
$x$ in $W_i $ can be identified with a neighborhood of a point of
$W_v$ so that
(locally) $S_{ji}$ is identified with $\MaxBar h_v$. In
particular $S_{ji}$
is also smooth and equidimensional.
\smallskip
Recall also that condition 5 in (1.13) says that, if the family is
j-compatible with the algorithm,
  $h_j(x) \leq h^{(t)}_j(x)$, and that
equality  holds if and only if $S_{ji}$ meets $(W_j)^{(t)}$ transversally
at $x$ (in 1.14 we explained the identifications involved in this
formula.)
\smallskip

(ii) Assume, as in (i) that (1.10.1) is j-compatible, and hence that
  $ \pi: W_j \to T $
is smooth. Recall that $ S_{jl} $ is a regular subscheme of $ W_j$ and
fix a point $ x \in S_{jl} $. Assume that
$ t \in T $, the image of $ x $, is a closed point let and $
(W_j)^{(t)} $ be the fiber over $ t $.
  These regular subschemes of $ W_j$ intersect at
$ x \in (W_j)^{(t)} \cap S_{jl} $, and it by 1.12(a) this
intersection is transversal at $ x $ if and only if  the restricted
morphism $ \rho_{jl}: S_{jl} \to T $, is smooth at $ x $.
\smallskip

(iii) Fix a smooth morphism $\pi : W \to T$, with $T$
integral, and a regular closed subscheme $C \subseteq W$ such that
the induced
morphism $q: C \to T$ is surjective, proper, and smooth. Then, for any $t
\in T$, the
number of
connected components of the geometric fiber $q^{-1}(\bar t)$ is constant.
\smallskip
This is probably known, but we should sketch a proof. The main steps are:
a) prove that all the fibers are irreducible, under the extra
assumptions:
  $C$ is irreducible and $q_{\ast}({\cal O}_C)= {\cal O}_T$
(That the fibers are all irreducible follows from the Stein
factorization theorem.)
b) Use a) to proof the assertion in case where any irreducible
component
$C_i$ of $C$ satisfies: $q_{\ast}({\cal O}_{C_i})= {\cal O}_T$.
c) The assertion reduces now to the following result, left to the reader:
Let $f:X \to Y$ be a smooth, proper morphism of integral
noetherian schemes, with Y normal.  Assume that $f_{\ast}({\cal
O}_{X}) \not= {\cal O}_Y$. Then there is an integral noetherian
scheme $Y'$ and a  finite, surjective morphism
$g: Y' \to Y$ such that if $X' \to Y'$ is the pull-back of $X$ via
$g$, we have
that $X'$ is not irreducible.
\medskip

(iv) Under the conditions of (iii) (i.e, $C$ is regular and $q: C \to T$
smooth), if we blow-up $W$ along $C$, to get
$g:W_1 \to W$, then $W_1$ is regular and $\pi_1:W_1 \to T$ is smooth. If,
in addition $\pi:W_0 \longrightarrow T $ is proper, the induced morphism
$\pi_1:W_1 \to T$ is smooth and proper and, if $t \in T$ is closed,
the fiber
${\pi _1}^{-1}(t)$ can be identified to the
blowing-up of
${\pi}^{-1}(t)$ with center $q^{-1}(t)$.
\smallskip

(v) Note that, given a family ${\cal G}$ as in (i), with $T$
integral, then
${\cal G}$ satisfies condition (AE) (resp. condition $\tau$) if and
only if, for all $t$ in $T$, the localization
${\cal G}\{t\}$
(see 1.10 (a)) satisfies condition (AE) (resp. condition $\tau$.)
First, let us check this assertion for condition (AE). Assume ${\cal
G}$ satisfies this condition. For each localization, the center
${C\{t\}}_i$
in the p-sequence of ${\cal G}\{t\}$ can be identified to the
$T\{t\}$-scheme obtained from $\rho _i:C_i \to T$ via the base
change $T\{t\}\to T$ (the natural morphism). Since $\rho _i$ is
smooth, proper and surjective, the same holds for
${\rho}\{t\} _i :{C}\{t\} _i \to {T}\{t\}$, i.e., ${\cal G}\{t\}$
satisfies (AE). For the converse, note that, since the center $C_i$
is regular, for all $i$ (1.13 (2)), the smoothness of $\rho _i$ is
equivalent to the fact that that for all $t \in T$ the fiber
${{\rho}_i}^{-1}(t)$ is a smooth $k(t)$-scheme, of the same
dimension. But since the fibers can be identified to fibers of the
localizations, this immediately follows from the smoothness of
${\rho}\{t\} _i$, for all $t \in T$.

The statement about condition $\tau$ follows from the fact that the
fiber of
${\cal G}$ at $t$ can be identified to the fiber of the
localization ${\cal G}\{t\}$, for all $t \in T$, and the general
fiber of ${\cal G}\{t\}$ can be identified to the general fiber of
${\cal G}$.
\smallskip

To finish this section, we mention some results about families of
embedded schemes.
\medskip

{\bf (2.6) Definition.}  If ${\cal F}=(X,W,\pi)$ is a family of
embedded schemes (1.10.3) we shall say that it is {\it equisolvable}
if the associated family of ideals (1.10.2) is equisolvable (2.4).
\smallskip

\proclaim
(2.7) Proposition. Assume ${\cal F}=(X,W,\pi)$ is an equisolvable
family of embedded schemes, where for all $t \in T$ the fiber $X_t$
is reduced . Then, for all $t\in T$, the resolution sequence (1.5.1)
induces the resolution sequence of $(X_t, W_t)$.

{\it Proof.} By definition, the equisolvability of ${\cal F}$ means
that the associated
family of ideals (see 1.10.3) ${\cal G}=(\pi:W \to T, I(X),
\emptyset)$ is equisolvable. Letting
$W_0=W$, ${\cal I}=I(X)$, $E_0=\emptyset$, consider the
p-sequence corresponding to the id-triple
${\cal T}=(W_0, {\cal I}, \emptyset)$. Equisolvability of ${\cal G}$
asserts that the p-sequence of ${\cal G}$ induces the p-sequence of the
fiber
${\cal G}(t):=(W_0^{(t)}, {\cal I}_0^{(t)}, E_0^{(t)})$,
in the following sense. We have
$C_0^{(t)}=C_0 \cap W_0^{(t)}$, $h_0^{(t)}=h_0|W_0^{(t)}$ and the
strict transform of $W_0^{(t)}$ in $W_1$ can be identified with
$W_1^{(t)}$ (see 1.9). Via this isomorphism,
$h_1^{(t)}=h_1|W_1^{(t)}$, $C_1^{(t)}=C_1 \cap W_1^{(t)}$, and so
on. Eventually, after $r$ steps, both the p-sequence of ${\cal T}$
and ${\cal G}(t)$ simultaneously stop. Now,
$I(X){\cal O}_{W_0^{(t)}}=I(X_0^{(t)})$ and, by the construction of
the associated algorithm for ${\cal S}$-couples (1.6), to obtain a
desingularization sequence for
$(X_0^{(t)},W_0^{(t)})$ we use the p-sequence of the fiber
${\cal G}(t)=(W_0^{(t)},{\cal I}_0^{(t)}, \emptyset)$, as
explained in 1.6. It is clear that, via the identification of
$W_i^{(t)}$ with a suitable subscheme  of $W_i$, for all possible
index $i$,
   $X_i^{(t)}$ corresponds to the strict transform of $X_0^{(t)}$ to
$W_i$. Moreover, the length $s_t$ (see the end of 1.6) of the
resolution sequence of
$(X_0^{(t)},W_0^{(t)})$ is constant, equal to the length $s$ of the
resolution sequence of $(X,W)$. In fact, $s$ (resp. $s_t$) is the
unique index such that the proper transform
$X_{s} \subset W_{s}$ (resp. $X_{s}^{(t)} \subset W_{s}^{(t)})$ has
the same codimension as the center
$C_{s}$ (resp. $C_{s}^{(t)}$). But
codim ($X_i^{(t)},W_i^{(t)})=$codim ($X_i,W_i)$ and, since $\cal G$
is equisolvable,
codim ($C_i^{(t)},W_i^{(t)})=$codim ($C_i,W_i)$. This proves our
contention about the indices. Since the desingularization functions
satisfy: $f_i=h_i$ (resp. $f_i^{(t)}=h_i^{(t)})$, $i=1, \ldots,
s=s_t$, the proposition is proved.
\bigskip

{\bf 3. PROOF OF THEOREM (2.3).}
\medskip

Now we begin the proof of  Theorem (2.3). For one implication, we show:
\proclaim
(3.1) Proposition. If our family ${\cal G}$ as in theorem~2.3
satisfies Condition (AE) then it satisfies
Condition $\tau$ and, moreover,  the p-sequence of
$  (W_0, {\cal I}_0, E_0)=(W, {\cal I}, E)$
induces that of any fiber ${\cal G}( t):=(W^{(t)}, {\cal I}^{(t)},
E^{(t)})$ or any geometric fiber.

By using 1.11(4) and 2.5(v), we see that it suffices to prove 3.1
for a local family, i.e., where $T= {\rm Spec}\,R$, with $R$ a
noetherian regular local ring. Concerning these, we can make the
following simplifying observations.
\medskip

{\bf (3.2) Remarks}. (a) Let $T= {\rm Spec}\,R$, with $R$ a
noetherian local ring. Using the fact that $T$ has a unique closed
point 0 (which belongs to any closed subset of $T$) plus other basic
results, it is easy to show the following fact.
Let $\phi :Z \to T$ be a proper morphism of schemes. Then,
  $\phi (Z)$ contains the closed point $0$ and if $\phi$ is smooth
at each point ${\phi}^{-1}(0)$, then $\phi$ is smooth everywhere on
$Z$.
\smallskip

(b) As an application, if ${\cal G}$ is a local family (1.10.1)
with p-sequence 1.2.1, and each induced morphism ${\rho}_i:C_i \to
T$ is proper, then: (i) for all $i$, $h_i$ reaches its maximum at
some closed point of the closed fiber of the induced projection $W_i
\to T$, (ii) if ${\rho}_i$ is smooth at each point of $C_i$ lying
over $0 \in T$, then ${\rho}_i$ is smooth.
\smallskip

(c) If ${\cal G}$ is a local family as in (b) (i.e., each $\rho _i$
is proper) and $z $ is the generic point of $T$, then 1.13(4) easily
implies that
$\Max{h_i}=\Max{h^{(z)}_i}$, where $h^{(z)}_i$ is the $i$-th
resolution function of the generic fiber
$(W^{(z)},{\cal I}^{(z)},{E}^{(z)})$.

Now we check a lemma. See (1.9) for an explanation of the
meaning of the equalities in {\it (i)} and {\it (iv)}.
\proclaim
(3.3) Lemma . Let ${\cal G}$ (1.10.1) be a local family which
satisfies the
general
hypotheses of Theorem 2.3, condition (AE) and is $j$-compatible, for
some $j$, $0 \leq j \leq r$ (where $r$ is the length of the
principalization sequence of (1.2.1). Let $t$ be
the closed point of $T$. Then :
\smallskip
(i) $\Max h_j = \Max h_j^{(t)}$ and,
if $x \in {\pi}_j ^{-1}(t)$, then
$h_j(x)=h^{(t)}_j(x)$.
\smallskip
(ii) For all $u \in T$, the number $c_j^{(u)}$ of irreducible
components of
the geometric fiber $C_j^{(\bar u)}$ is constant.
\smallskip
(iii) $W_j^{(t)}$ (identified to
$\pi _j^{-1}(t)$) intersects $C_j$ transversally,
\smallskip
(iv) ${\rho}_j^{-1}(t) = C_j \cap W_j^{(t)} = C_j^{(t)}$ (scheme
theoretic
intersection).
\smallskip
(v) If $\pi _{j+1}: W_{j+1}\to T$ is the induced morphism, then
$\pi _{j+1}$ is smooth and
${\pi _{j+1}}^{-1}(t)$ is canonically isomorphic to the blowing-up
of $W_j^{(t)}$ along $C_j^{(t)}$.
\medskip

{\it Proof.} (i) Fix $ j $ as above and let
$ a_{j,0}, a_{j,1},\ldots, a_{j,s}$ be all the possible values of
$h_j: W_j \to I^{(n)}$ ($n=$ dim $W$.) Let
$S_{jl}:=f_j^{-1}(a_{jl}),\, i=0, \ldots, s$. This is, with the
reduced structure, a locally closed and regular subscheme of $W_j$
(1.13(2)). First we
shall see that, via the naturally induced projection, $S_{jl}$ is
smooth over $T$, for all $l$. We argue as in 2.5 (i), but now
we are
concerned with smoothness over $ T $. It is enough to check this
smoothness
  on a
neighborhood (in $W_j$) of an arbitrary point
$x \in S_{jl}$ lying over $t$. So, take such a point.
Now, by property (iv) of 1.2,
there is a largest index
$v \in \{j, \ldots, r \}$ so that a neighborhood of $ x $ in $W$
can be identified
with an open subset of $ W_v $; this identification is such that
  $x$ corresponds to a point
$x'$ of $C_v$ and $ S_{jl}$ is (locally) identified with $C_v$, moreover
this identification is compatible with the projection on $ T $.
By condition (AE) $C_v$ is smooth over $T$, so we
have that
$S_{jl} \to T$ is also smooth, proving our claim. From this,
(2.5)(ii) insures that $W_j ^{(t)}$ and $S_{jl}$ meet transversally
at any common point, for
$l \in \{j, \ldots, r\}$, and hence by (1.13)(5),
$$ (3.3.1)\qquad h_j ^{(t)}(x)=h_j(x), {\rm for\, any }\, x \in
(W_j)^{(t)},$$
showing the second part of (i).
To see the first part, let $t
\in T$ and $x \in {\pi}_j ^{-1} (t)$ (identified with $W_j ^{(t)}$),
be such that
$h_j ^{(t)}(x)=a_j ^{(t)}:=\Max h_j ^{(t)}$. Then, by (3.3.1) and
Remark (2.5) (i), we have
  $a_j ^{(t)} := h^{(t)}_j(x)=h_j (x) \leq  a_j:=\Max h_j  $. Now let
$w \in {C_j \cap W^{(t)}_j}$.
Since $C_j$ is the set of
points of $W_j$ where $h_j$ reaches its maximum value $a_j$, we have
$a_j=h_j(w) = h^{(t)}_j(w)\leq a^{(t)}_j$. The desired equality is
proved.
\smallskip

(ii) Since we assume $C_j$ proper over $T$ we may use (2.5) (iii) to
conclude that $c^{(u)}_j$ is constant, for all $u$ in $T$.
\smallskip

(iii) This is (3.3.1) together with (2.5) (ii).
\smallskip

(iv) This is a consequence of (3.3.1) and the transversality proved
in (iii).
\smallskip

(v) The first part follows from  Remark 3.2, the second is 2.5(iv).
\medskip

{\bf (3.4).} Now we prove Proposition (3.1). As already seen we may
assume that ${\cal G}$ is local, let $t$ again denote the closed
point of $T$. Let $r$ be
  the length of the principalization sequence of $(W, {\cal I}, E)$. We
proceed step-wise. Since every family is 0-compatible with the
algorithm, we may apply Lemma (3.3) to our ${\cal G}$, with $j=0$.
We get:
\smallskip

($\alpha _0$) $\Max h_0 ^{(t)}=\Max h_0$ and if $\pi _0(x) = t$ then
$h_0 ^{(t)}(x)= h_0(x)$,
\smallskip

$(\beta _0)$ $c_0 ^{(u)}$ is constant, for $u \in T$,
\smallskip

($\gamma _0$) $C_0 \cap W_0^{(t)} = C_0^{(t)}$,
\smallskip

($\delta _0$) $W_1^{(t)}$ gets identified to
$\pi_1^{-1}(t) \subset W_1$.
\smallskip

By assumption, $C_0 \to T$ is smooth
and onto; this together with ($\gamma _0$) show that ${\cal G}$ is
also 1-compatible. If $r \geq 1$, apply Lemma (3.2) again, with
$j=1$. We get:
\smallskip

($\alpha _1$) $\Max h_1 ^{(t)}=\Max h_1$,if $\pi _1(x) = t$ then
$h_1 ^{(t)}(x)= h_1(x)$,
\smallskip

$(\beta _1)$ $c_1 ^{(u)}$ is constant,
\smallskip

($\gamma _1$) $C_1 \cap W_1^{(t)} = C_1^{(t)}$,
\smallskip

($\delta _1$) $W_2^{(t)}$ gets
identified to
$\pi_2^{-1}(t) \subset W_2$.
\smallskip

Again we get, as above, that
${\cal G}$ is 2-compatible. If $r \geq 2$, repeat. Eventually we
get, for all $u \in T, i \in \{0, \ldots,r \}$:
\smallskip

($\alpha$) $\Max h_i ^{(t)}=\Max h_i$,
\smallskip

($\beta$) $c_i ^{(u)}$ is constant,
\smallskip

($\gamma$) $C_i \cap W_i^{(t)} = C_i^{(t)}$,
\smallskip

($\delta$)  $W_i^{(t)}$ gets identified to
$\pi_i^{-1}(t) \subset W_i$ and, using this identification,
if $\pi_i(x) = t$ then
$h_i ^{(t)}(x)= h_i(x)$
\smallskip

Since $h_r$ is constant along $W_r$, it follows that
$h_r^{(t)}$ is constant along $W_r^{(t)} \subset W_r$, i.e.,
$r_t=r$. Noticing that 1.13(4) implies that if $z$ is the generic
point of $T$ then
$\Max h_i=\Max h_i^{(z)}$ (3.2(c)), now it becomes clear that
$(\alpha)$ and
($\beta$) say that that condition $\tau$ holds for $\cal G$, and
($\delta$) says that the principalization sequence of
$  (W_0, {\cal I}_0, E_0)=(W, {\cal I}, E)$
induces that of the closed fiber ${\cal G}( t):=(W^{(t)}, {\cal I}^{(t)},
E^{(t)})$.
The corresponding statement for geometric fibers
follows from (1.13) (1). This proves Proposition 3.1.
\medskip

Next we shall prove:
\proclaim
(3.5) Proposition. If a family ${\cal G}$ (as in Theorem 2.3) satisfies
condition
$\tau$, then it also satisfies (AE).

First, we shall prove:

\proclaim
(3.6) Lemma.
Let ${\cal G}=(\pi:W \to T, {\cal I}, E)$ be a local family, with
$T$ integral (where $t$ is the closed point of $T$),  satisfying
Condition $\tau$. Let 1.2.1 be the principalization sequence of $(W,
{\cal I}, E)$,
where the projection $\rho _i:C_i \to T$ is proper
for all $i$. Moreover, we assume
that the family is $j$-compatible with the algorithm, for some $j
\in \{0,
\ldots, r \}$.
Then:
\smallskip
(i) $\Max h_j = \Max h_j^{(t)}$,
\smallskip
(ii) $\rho _j:C_j\to T$ is smooth and surjective,
\smallskip
(iii) $ C_j \cap W_j^{(t)} = C_j^{(t)}$, and this is transversal,
\smallskip
(iv)  If $ W_{j+1} \to W$ is the blowing-up of $W$ along $C_j$ and
$\pi_{j+1}: W_{j+1}\to T$ the induced morphism,
then $\pi _{j+1}$ is smooth and the
subscheme ${\pi _{j+1}}^{-1}(t)\subset W_{j+1}$ is canonically
isomorphic to the blowing-up
of $W_j^{(t)}$
along $C_j^{(t)}$.

{\it Proof.} (i) Let $a=\Max h_j$. Then (see 3.2(c)) if $z$ is the
generic point of $T$, $\Max h_j=\Max h_j^{(z)}$. But by condition
$\tau$, $\Max h_j^{(t)}=\Max h_j^{(z)}=a$.
\smallskip

(ii) By 3.2(b), to show that $\rho _j$ is smooth it suffices to
show that $\rho _j$ is smooth at each point $x$ lying over $t$, and
by 2.5(ii) this happens if and only if $h_j(x)=h_j ^{(t)}(x)$. Now,
by 1.13 (5), always
$h_j(x) \leq h_j ^{(t)}(x)$. On the other hand, since $x \in C_j$,
$h_j(x)=\Max h_j=a$. We have
$h_j ^{(t)}(x) \leq \Max h_j ^{(t)}(x)=a=h_j(x)$, where we have used
(i). Hence,
$h_j(x)=h_j^{(t)}$. By
3.2(a) the morphism $\rho _j$ is onto, so (ii) is proved.
\smallskip

The fact that the intersection is transversal was seen in (ii). In
particular, the scheme-theoretic intersection (i.e., what is meant
in (iii)) is a reduced, regular subscheme of $W_j$. By (i), the
inclusion
$C_j \cap W_j ^{(t)} \subseteq C_j ^{(t)}$ is valid. We shall see
that a strict inclusion is impossible. Assume, by contradiction,
that $C_j \cap W_j ^{(t)} \subset C_j ^{(t)}$.

By (2.5) (iii), all the geometric fibers
of $C _j : \to T$
(which is smooth, by (ii), and proper) have the same number of
irreducible
components, say $c'$. On the other hand, for all $t \in T$, by Condition
$\tau$ the number of components of
$C_j^{(t)}$ must be the same, say $c$.
Now, by using the Generic Smoothness Theorem ([14], p. 272)
(applied to the different projections
$S_{i j} \to T$, $S_{i j}= h_j ^{-1}(a_{ij})$ (where
$a_{j1} < \cdots <a_{js}:=a$
are all the values of $h_j:W_j \to I^{(n)}$), we see
that for a non-empty open set $U \subseteq T$, the restriction of
the family to $U$ satisfies condition (AE). The generic point $z$
of $T$ is in $U$. Hence, by Lemma 3.3,
$C_j \cap W_j ^{(z)} = C_j ^{(z)}$, so $c=c'$. But by 1.13(2), for
dimension reasons, $C_j \cap W_j ^{(t)}$ must be a union of
connected (or irreducible) components of $C_j ^{(t)}$. Since we are
assuming
$C_j \cap W_j ^{(t)} \subset C_j ^{(t)}$, we get $c < c'$, a
contradiction.
\smallskip

(iv) This is seen by using (2.5)(iv).
\medskip

{\bf (3.7)} Now we prove Proposition (3.4).
By (2.5 (v)), it suffices to show each localization
${\cal G}\{t\}$, $t \in T$,  satisfies condition (AE). So we may
assume ${\cal G}$ local (with $t$ the closed point of $T$.)
So,  Let ${\cal G}=(\pi:W \to T, {\cal I}, E)$ be as in  (1.10.1)
and local, a
let $r$ be
  the length of the principalization sequence of
$(W , {\cal I}, E)$.

Again, we proceed step-wise. Since every family is 0-compatible with
the algorithm, we may apply Lemma (3.6) to our family ${\cal G}$,
with $j=0$. We get, for all $t \in T$:
\smallskip

($\alpha _0$) $C_0 \to T$ is smooth and surjective (and, by
assumption, proper),
\smallskip

($\beta _0$) $C_0 \cap W_0^{(t)}$ is transversal.
Thus the family is also 1-compatible.
\smallskip

If $r \geq 1$, we apply (3.4)
again, to get:
\smallskip

($\alpha _1$) $C_1 \to T$ is smooth, surjective and, proper,
\smallskip

($\beta _1$) $C_1 \cap W_1^{(t)}$ is transversal.
\smallskip

If $r \geq 2$, repeat. Eventually we get that the induced projection
$C_i \to T$ is smooth, proper and surjective for all
$i \in \{0, \ldots , r \}$, i.e., the family ${\cal G}$ satisfies
Condition (AE). Proposition 3.5 is proved.
\smallskip

Clearly, Theorem (2.3) follows from propositions (3.1)
and (3.5).
\medskip
{\bf (3.8) Remark.}
Recall that the family $ \cal G $ is {\it
equisolvable} if
the equivalent conditions hold (2.4). But there is more to say about this
notion which follows from the proof of the equivalence:
\smallskip
{\bf 1.} If the equivalent conditions of Theorem 2.3 hold, then
$h_j(x) = h_j^{(t)}(x)$, ($x \in W_j^{(t)}$), for any
$t \in T$ and any index $ j=0,1,...,r $ in the principalization sequence
(1.2.1).

This follows from the proof of iii) in 3.6, noting that if the equivalent
conditions hold, then we may take $ U=T$ and argue as we did there,
but now
for any index $j$.
\smallskip
{\bf 2.} Suppose that
$ \cal G $  is a equisolvable family, let $ T_1 \to T $ be a morphism of
regular schemes and $ {\cal G}_1 $ be the induced family over
$ T_1 $. Then: a) $ {\cal G}_1 $ is equisolvable, and b) the
principalization  sequence of $ ( W_1, { \cal I}_1, E_1) $ is the
pull-back
of the principalization  sequence of
$(W,{\cal I},E)$.
\smallskip
In fact, a) follows since condition $ \tau $ is a condition on the
fibers,
while b) follows from the proof of (3.6) (iii) which shows that the
principalization sequence of $ ( W, { \cal I}, E) $ is obtained by
``putting
together'' the principalization sequences of the different fibers $ (
W^{(t)}, { \cal I}^{(t)}, E^{(t)}) $ .
\smallskip
So what we have achieved, via the equivalence in Theorem 2.3, is a
condition of equisolvability of a given family, which can be expressed
entirely in terms of the fiber. This will be a key point in the
applications to be discussed in the next chapter.
\bigskip

{\bf 4.  STRATIFICATIONS OF HILBERT SCHEMES.}
\medskip

In this section we intend to show that given a family of ideals
(1.10.1) it is possible to naturally express the parameter scheme $T$
as  a disjoint union of locally closed sets, such that the
restriction of the family to each of these is equisolvable. Because
of definition  2.6, this immediately implies a similar result for
families of embedded schemes (1.10.2). An interesting application is
to the universal family parametrizing subschemes of a fixed projective
variety  (by the theory of Hilbert schemes). To accomplish this we
need a slight generalization of the notion of family introduced in
1.8 (see 4.7).

Throughout this section we retain the assumptions of 2.1, i.e., we
work within a given allowable collection ${\cal S}$ where a good
principalization algorithm has been defined.
\medskip

{\bf (4.1) Definition} . A function
$\gamma$
from an algebraic variety $T $ to a partially ordered set $ (I,
\leq )$ is
called an
  {\it LC-function} if it has the following property: for any irreducible
closed subset
$ Z $ of $ T $, there is a dense Zariski open subset $ U $ of $ Z $
such that
\smallskip

(a)  $\gamma$ is constant along points of U ( let $\gamma_1 $ denote
such a constant value),
\smallskip

(b) at any point $ x \in Z $: $ \gamma_1 \leq \gamma(x) $.
\smallskip

In our context $ (I, \leq )$ will be totally ordered and one can check
directly that a) and b) hold iff  $\gamma$ is upper-semi-continuous and
takes only finitely many values.
This notion appears in
[24], page 241. It can be proved that  the fibers of such a $\gamma$
define a
  finite partition of $Z$ into locally closed subsets (Theorem 2.1 in
[24]).
  Recall from (2.2) that, given a family ${\cal G}$ (1.10.1), we have
defined a
  function $\tau _{{\cal G}}$ from the parameter space
$T$ into the totally ordered set ${\Lambda ^{(m)}}$, where $m$
denotes the
dimension of the fibers of $W \to T$.
\smallskip

Now we have the following:

\proclaim
(4.2) Theorem . Consider a family
$ {\cal G}=(\pi:W \to T, {\cal I}, E)$
as in (1.10.1) such that all
morphisms
$\rho_i: C_i \to T $ are closed ([14], p. 100) , $i=0,\ldots,r$
(2.1). The function
$\tau := \tau _{{\cal G}}$ from $T$ to ${\Lambda ^{(m)}}$ is an
LC-function, i.e.
(a) and (b) of (4.1) hold.

{\it Proof of part (a).} First, some reductions valid for both
parts. Let $ Z $ be an irreducible closed subset of $ T $.
Since the function $ \tau $ is defined entirely in terms of the fibers,
we may replace $ Z \subset T $ by an embedded
desingularization  $ Z_1 \subset T_1 $ and
the family ${\cal G}$
by its pull-back via the induced morphism $Z_1 \to T$. Hence, there
is no loss of generality to assume that, in 1.10.1, our closed set
$Z$ is the whole regular scheme $T$, and that $T$ is irreducible.

Now we check a) Let (1.2.1) be the principalization sequence of
$(W, {\cal I}, E)$   .
For each index
$ i=0,\ldots,r $ the image of $ C_i $ in $ T $ is a closed subset.
Replacing $ T $ by an open subset (still denoted by $ T $), we may
assume that
all $ C_i $ map surjectively on $ T $.
Each $ C_i $ is a regular subscheme in $ W_i $. Since the morphism from
$ C_i $ to $ T $ is of finite type, and the field of rational
functions of $ T $
has characteristic zero, it follows from the Generic Smoothness
Theorem that the morphism is smooth
over the generic point of $ T $. Hence there is an open dense set $
U \subset T $ over
which the morphism is smooth. So ultimately, after suitable
restriction of
$ T $ we may assume that each $ C_i $ maps properly and smoothly on
$ T $.
Now a) follows from Theorem 2.3.
\smallskip
The proof of property b) will be given in (4.6) after some
previous results are discussed.

\proclaim
(4.3) Lemma.  Let ${\cal G}$ (1.10.1) be a family of ideals, having
1.2.1 as p-sequence, $j \in \{1, \ldots, r \}$ an index, such that:
\smallskip
(1) $ T $ is regular, irreducible and $ \rho_s:C_s \to T $ is smooth and
surjective, for $s=0,1, \ldots,j-1 $.
\smallskip
(2) for $ s \in \{0,1,\ldots ,j-1 \} $, $ \Max h_s^{(t)}= \Max h_s$,
$ c_s^{(t)}$ is constant, and  $C_s^{(t)}=C_s \cap W_s^{(t)}$ ($
s=0,1,...,j-1$), for all $ t \in T $ (where we have used
the convention of the Remark 1.12).
\smallskip
(3)  $\rho_j: C_j \to T $ is surjective.
\smallskip
Then:
\smallskip
(i) $ \Max h_j^{(t)} \geq \Max h_j$, for any $t \in T $,
\smallskip
(ii) if equality holds in (i) for a certain $ t $, then the fiber $
\rho_j^{-1}(t)$ and the geometric fiber
$ \rho_j^{-1}(\bar{t})$ are regular, and $ c_j^{(t)}\geq
\alpha_{j,t}$, the
number of connected components of $ \rho_j^{-1}(\bar{t})$.

{\it Proof.} First note that assumption (1) implies the smoothness of the
projection $ W_j \to T $. Now let $t \in T $ and $ D^{(t)}_j=
\rho^{-1}_j(t)$. By (3), $ D^{(t)}_j \neq \emptyset $. At any point
$x$ of
$ D^{(t)}_j $, by 1.13 (5),
$ h_j (x) \leq h^{(t)}_j(x)$, with equality if and only if
$C_j$ and
$ W^{(t)}_j$ meet transversally at $x$, which (by 2.5 (i) and (ii)) is
equivalent to the smoothness of $ \rho_j: C_j \to T $ on a
neighborhood of
$x$. (Here we are using the conventions of 1.9 and 1.13 (5).) In
particular, since $h_j(x)= \Max h_j $ for any $ x \in C_j$, we
obtain (i). Moreover, if equality holds, then
$$ \Max h_j=h_j(x) \leq h^{(t)}_j(x) \leq \Max h^{(t)}_j = \Max h_j $$
hence $ h^{(t)}_j(x)=h_j(x)$ so $ \rho_j: C_j \to T $ is smooth near $x$.
Thus $ \rho_j $ is locally smooth in this case, and we get an
open neighborhood $ V $ of $ D^{(t)}_j $ such that $ \rho_j$  is
smooth on $V$.

Now, $ D^{(t)}_j  \subseteq C^{(t)}_j $. It also follows, as in the
proof of
3.6 (iii), that both schemes are equidimensional and of the same
dimension.
Hence , $ D^{(t)}_j $ is a union of connected components of $
C^{(t)}_j $. Thus
$ D^{(t)}_j $ is regular (because  $ C^{(t)}_j $ is so), and
  $ c^{(t)}_j \geq \alpha_{j,t}$. This  proves (ii).
\medskip

{\bf (4.4)} In the sequel, we shall use the following notation.
Given a family ${\cal G}=(\pi:W \to T, {\cal I}, E)$ (1.10.1),
$z \in T $, write
$$ \tau_{\cal G}(z)=(a_0 (z), b_0 (z),
a_1 (z), b_1 (z), \ldots,a_{r(z)} (z), b_{r(z)}
(z),\infty,\ldots,\infty,\ldots) $$
where $a_i (z)= \Max h_i^{(z)} $, $b_i (z)= c_i^{(z)} $, for all i.

Let $ U
$ be a dense open set of $ T $ where the restriction of $ {\cal G}$
to $ U
$ satisfies condition AE ( or $ \tau $), as given in the proof of
4.2 (a).
Note that for $ t \in U $, $\tau_{\cal G}(t)=\tau_{\cal G}(t_g) $,
where $
t_g $ is the generic point of $ T $.
\medskip

\proclaim
  (4.5) Lemma.  Let ${\cal G} $ (1.10.1) be a family of ideals with
p-sequence 1.2.1, where we fix an index $j$ and assume conditions (1) and
(2) of (4.3). Suppose now that  $ \rho_j: C_j \to T $ is
proper but not surjective,
let
$ z \in {\rm Im} (\rho_j) $. Then $ a_j (z)>a_j (t)$ for $ t \in U
$ (where $U$ is as above in (4.4)).

{\it Proof.} Consider an open dense set $ U $ in $ T $ as above,
i.e. such that the restriction ${\cal G}'$ of ${\cal G}$ to $U$
satisfies condition (AE). Then, certainly $ z \notin U$. To determine
the value of $ a_j(t) $, let $ j' $ be the smallest index
such that $ j < j' \leq r $  and $ \rho_{j'}: C_{j'} \to T $ is
surjective. Such
$ j' $ exists, since $ C_r = W_r $ (because the last
principalization function $h_r$ is constant on $W_r$, cf.  1.2 (ii))
and the projection $ W_r
\to T $ is onto. It follows from the definitions that the step
number $ j $ in the
principalization
sequence of ${\cal G'} $ is the restriction of the step number $ j' $ of
the principalization sequence of ${\cal G} $. From this we easily get:
$$ a_j(t)= \Max h_{j'}$$
>From 3.3(i) we get $h_s(x)=h_s^{(t)}(x)$ if $x \in W_s$ lies over
$t \in U$, $0 \leq s \leq r$.
Now choose $ x \in C_{j'} $ mapping to $ z $. Then we have
$ h_{j'}(x)=\Max h_{j'}$, and
$$ a_j (z):=\Max h^{(z)}_j > \Max h^{(z)}_{j'} \geq
h^{(z)}_{j'}(x) \geq h_{j'}(x) = \Max h_{j'}= a_j(t) $$
where we have used 1.2 (v) and 1.13 (5). This shows the desired
inequality.
\medskip

{\bf (4.6)} {\it Proof of (4.2), b)}. As indicated in (4.2), we may
assume that
$ Z= T $ is regular and irreducible. We use the notation of (4.4).
Moreover, let
$ 0 \leq i_0 <i_1 <i_2 < \ldots <i_r =r $ be the indices such that the
corresponding projection $ \rho_{i_j}: C_{i_j} \to T $ is surjective
(dominant).
We assume that all $ \rho_i $ are proper.

Recall our definition of the open dense set $ U $ in $ T $ (cf. proof of
4.5 or 4.2 a)) and note that,

1) $ \rho_s(C_s)\cap U \neq \emptyset $ iff
  $ s \in \{ i_0,i_1, \ldots,i_r \}$

2) the index $i_j$ induces the index $j$ over $U$, and
$\Max h^{(t)}_j =\Max h_{i_j} $ for $ t \in U $ (see (1) in Remark 3.7).

First, we compare $ a_0 $ at a point $ z \in T $ and at a point $ t
\in U $.
Let $ s $ be the
smallest index so that $ z \in {\rm Im} (\rho_s) $. Then:

(i) $s\leq i_0 $,

(ii) $ a_0 (z):=\Max h^{(z)}_0 \geq \Max h_s $ ( choose $ x \in C_s $
mapping
to $ z $ and use properties (iv) in 1.2 and (4) in (1.11) to show that
$ h^{(z)}_0
(x) \geq h_s (x)=\Max h_s $),

(iii) $ a_0 (t)=\Max h_{i_0} $

If $ s < i_0 $ then $ a_0(z)>a_0(t) $ by ii) and iii), together with
(1.2),(v). Hence $ \tau_{\cal G}(z)>\tau_{\cal G}(t)$ and (b) is
proved in
this case.

If $ s= i_0 $, we may replace $ T $ by a suitable open neighborhood
$ T' $
of $ z $ so that $ s= i_0= 0 $; in fact set $ F_e $ as the union of all
closed subsets
$ {\rm Im} \rho_k $, for $ 0\leq k < e $, hence $ z $ is not in $
F_{i_0} $. Set
$ T'= T- F_{i_0} $.
In this case we may apply Lemma (4.3) for $ j=0 $, in fact all
requirements in (4.3) are
trivially fulfilled for $ j=0 $.
Hence $\Max h^{(z)}_0 \geq\Max h_0 $
by (4.3) (i). On the other
hand, applying 2) (above) for $ t \in U \cap T' $ it follows that
$\Max h^{(t)}_0 =\Max h_{0} $ and hence that $ a_0(z) \geq
a_0(t)$. If the
inequality is strict then $ \tau_{\cal G}(z)>\tau_{\cal G}(t)$ for any
$ t \in U \cap T' $. If not, then $\Max h^{(z)}_0 =\Max h_{0} $,
hence (ii) in Lemma (4.3) applies so
$ c^{(z)}_0 \geq \alpha_{0,z}$. By 2.5 (iii) the geometric fibers
have the
same number of irreducible components ( equal to $\alpha_{0,z}$) for all
$ t \in U \cap T' $. Hence, $ b_0(z) \geq b_0(t)$ for $ t \in U
\cap T' $.
If $ b_0(z) > b_0(t)$, we are done ($ \tau_{\cal G}(z)>\tau_{\cal
G}(t)$ ).
Otherwise $ (a_0(z),b_0(z))=(a_0(t),b_0(t))$
; we then
compare $ a_1(z) $ and
$ a_1(t) $. Let $s_1 $ be the smallest index $ j > s=i_0 $ so that
$ z \in {\rm Im}\,\rho_j $. Then $ s_1 \leq i_1 $.
If $ s_1 < i_1 $ we argue as in the case $ s < i_0 $ above to show that
$ a_1(z)>a_1(t) $. Hence $ \tau_{\cal G}(z)>\tau_{\cal G}(t) $ and we are
finished.
If $ s_1 = i_1 $, we may replace $ T $ by a suitable open
neighborhood $ T'
$ of $ z $ so that $ s_1= i_1= 1 $, in fact set $ F_{i_1} $ as the
union of all
closed subsets $ {\rm Im} (\rho_k) $, for $ 0 \leq k < i_1 $ and
$ k \neq i_0 $,
hence $ z $ is not in $ F_{i_1} $. Set $ T'= T- F_{i_1} $. Note now
that the
hypotheses  of Lemma (4.3) hold for $ j=1 $. So, by applying Lemma (4.3)
we get $
a_1(z) \geq a_1(t) $, with equalities implying $ b_1(z) \geq b_1(t)$ and
the smoothness of
$ \rho_{i_1}$ (restricted to a suitable neighborhood of $ z $). Clearly,
after repeating this process  at most $ r $ times, we obtain
$ \tau_{\cal G}(z)\geq \tau_{\cal G}(t) $, in the lexicographic order, as
desired.
\smallskip
In the definition of families of ideals  and of embedded
schemes (1.10, (a) and (b) respectively), we assumed that the
parameter space $T$ was a regular scheme. This is too restrictive
for what we plan to do next, so we introduce the following
definitions.
\medskip

{\bf (4.7) Definition.} (a) A {\it  general family of embedded
schemes} is
a pair
$$ (4.7.1) \qquad {\cal F} = (j: X \to W, \pi:W \to T) $$
as in (1.10.2), where $ \pi:W \to T $ is a smooth morphism of
equidimensional schemes over a field $k$ of characteristic zero and
$p:= \pi j:X
\to T$, is flat (but we do not
assume that $W $ or $ T $ are regular).
(all these
are schemes over a field $k$ of characteristic zero).

As usual
$X_t:=p^{-1}(t), W_t={\pi}^{-1}(t)$ denote the fiber over $t \in
T$. Note that
that dim $  W_t $ is constant, say $=m$, for $ t \in T $.
\smallskip
(b) A {\it  general family of ideals} is
a pair
$$(4.7.2)\qquad {\cal G}=(\pi:W \to T, {\cal I}, E)$$
where $ \pi:W \to T $ is a smooth as in (1.10.1) and all the
fibers
${\cal G}(t):=(W^{(t)}, {\cal I}^{(t)}, E^{(t)})$
are idealistic triples in ${\cal S}$ (1.1) (but we do not
require $ W $ or
$ T $ to be regular or irreducible).

Note that
$\tau _{{\cal G}}: T \to \Lambda^{(m)}$
still can be defined for any general family.
If $ T $ admits a desingularization $ T_1 \to T $
then $ {\cal G}_1 $
(obtained from ${\cal G}$ by base change to $T_1$) is a family in the
sense of definition
(1.10)(a). Hence  Theorem 4.2
applies to $ {\cal G}_1 $. Since  $\tau _{{\cal G}}$ is defined in terms
of the fibers only, from the properness of $ T_1 \to T $ we see that
the conclusion of (4.2) is also valid for the family $ {\cal G}$.
This shows:

\proclaim
(4.8) Theorem.
If ${\cal G}$ is a general family of ideals, then $\tau _{\cal G}$ is an
LC-function along T.

{\bf (4.9)} As in (1.10.c), we associate to a general family of
embedded schemes ${\cal F} = (j: X \to W, \pi:W \to T) $ (4.7.1),
a general family of ideals, say
${\cal G}=(\pi:W \to T, {\cal I}, E)$, where ${\cal I}=I(X)$ and $
E=\emptyset$. Using this notation,
  For simplicity, in the sequel we denote by $\tau _{{\cal F}}$ the
function $\tau _{{\cal G}}$ corresponding (see  2.2) to thefamily
of ideals
${\cal G}$.
\smallskip

Now fix the following objects:

(i) a graded algebra $S$ over a field $k$ of characteristic zero,
finitely generated by elements of degree one,
such that
$ W=Proj( S ) $ is smooth over $ k $, of dimension $n$, and it belongs
to the allowable class ${\cal S}$ of our algorithm;

(ii) an element $\alpha\in \Lambda^{(n)}$ and

(iii) a polynomial $ Q $ with rational
coefficients.

If T is a scheme of finite type over $ k $, set $ W_{T}:=W \times T $.
Let $ {\cal H}(W, \alpha, Q) $ be the
class of all
general families of reduced schemes
  ${\cal F}_{T}$, of the form $( j : X_T \to W_T, \pi) $ where
$ \pi: W_T \to T $ is obtained from $ W \to Spec(k) $ by base change,
and $ X_T \subset W_T $ is a closed subscheme such that the induced
projection
$ X_T \to T $ is flat ($j $ being the inclusion), so that:
\smallskip

(0) T is a reduced.
\smallskip

(1) For all $ t \in T $, $\tau _{{\cal F},T}(t)= \alpha $ (hence
the pull-back of this
general family via $T' \to T$, with $T'  \in {\cal S}$ regular is
{\it equisolvable} since it
satisfies condition $ \tau $.)
\smallskip

(2) If $ ( X_{t}, W) $ is the couple induced by $ ( X_{T}, W_{T}) $  over
$ t \in T $ then $ Q $ is the Hilbert polynomial (relative to the
line bundle corresponding to $S(-1)$) of the
embedded scheme $ X_{t} \subset W $.
\smallskip

Now we can state the following:
\proclaim
{\bf (4.10) Theorem}. Under the conditions of 4.9 (and letting
$ H:={\cal H}(W, \alpha, Q)$) there is a universal object in the
class of general families in $H$ (to be called an universal
$(\alpha,Q)$-equisolvable
family).
That is, there is a general family,
${\cal F}_{H(\alpha,Q)}$ defined by
$ ( X_{H(\alpha, Q)}\subset W_{H(\alpha, Q)}) $, $W_{H(\alpha, Q)}
\to H(\alpha, Q)$)
such that for each general family ${\cal F}_{T}$ in
$ H $, there is a unique morphism
$ T \to H(\alpha, Q)$ so that ${\cal F}_{T}$ is the pull-back of
${\cal F}_{H(\alpha ,Q)}$.

{\it Proof.} Consider the Hilbert scheme $ H(Q) $, parametrizing
subschemes of
the projective variety $ W $ having Hilbert polynomial $ Q $.
We refer here to [22]
(see (c) on page 21) for a summary of results on Hilbert schemes.
Let $ X(Q) \subset W \times H(Q)$ be the universal family; note that
$ (X(Q), W \times H(Q))$ together with the projection $W \times H(Q) \to
H(Q)$ defines a general family, say ${\cal F}(W,Q)$. Theorem~4.8
says that
$\tau _{{\cal F}(W,Q)}: H(Q) \to \Lambda^{(n)}$ is an LC-function. Hence
its fibers define a partition of $ H(Q) $ into a disjoint union of
locally closed
subsets.
If $\alpha \in \Lambda^{(n)}$, let
$ H(\alpha, Q):= [\tau _{{\cal F}(W,Q)}]^{-1}(\alpha)$.
Given a general family, say ${\cal F}_T $,  in
$ H $, since the natural morphism $ T \to H(Q) $ obtained by
universality of
Hilbert schemes has constant value $\alpha$, it becomes clear that
it factors though $ H(\alpha, Q) $ and vice versa,
proving the
theorem.
\medskip

{\bf (4.11) Remark.} We define an {\it LC-partition} of a scheme
$T$ as an expression of $T$ as a disjoint union of locally closed
subsets (or subschemes, with the reduced structure). Thus, in 4.10
we have introduced in the Hilbert scheme $H(Q)$ an LC-partition
naturally related to equiresolution, by means of the (reduced)
subschemes $ H(\alpha, Q) $. We might call this a  {\it
stratification} of $H(Q)$ but, in the literature, this term is often
used in a more restricted sense. For instance, in
  [20] the notion of {\it stratification} of
$T$ is introduced  as an LC-partition (each set thereof is called a
stratum)
with the properties:

(i) the boundary of each stratum is a union of strata;

(ii) the singular locus of the closure of each stratum is a union of
stratum,
and

(iii) each stratum is smooth.
\smallskip

In 2.5(c) of [20], it is remarked that given an
LC-partition of an algebraic variety T, one can attach to this partition
a {\it coarsest} stratification of $T$, with the property  each
locally closed
subset of the partition is a union of strata. Here
{\it coarsest} means that any other stratification with this
property is a
refinement of the first.
\smallskip
So, if we want to obtain a stratification (in the sense of Lipman)
of the Hilbert scheme $H(Q)$,  naturally related to equiresolution,
all what we need it to take the coarset stratification associated to
the LC-partition described in theorem 4.10.
\bigskip

{\bf 5. DESCRIPTION OF AN ALGORITHM OF PRINCIPALIZATION.}
\medskip

{\bf (5.1)} This chapter will be devoted to
recall some results on the desingularization process developed in
[11] and, for the
reader's convenience, to review, without
proofs, the basic ideas behind the process, as well as notation and
terminology which are essential in the discussion of section~6.
\medskip

As already mentioned (see 1.6), it seems convenient to initially
investigate
not the problem of embedded resolution of schemes or
principalization of an ideal, but rather a seemingly more technical
one, involving basic objects.

Namely, a {\it basic object} is a
triple $B=(W,(J,b),E)$ where $W$ and $E$ are as in 1.1, $J$ is a sheaf
of ideals of $\O_{W}$ such that $J_{\xi}\neq 0$ for any point $\xi\in
W$, and $b$ is a non-negative integer.
Let $\Sing(B)=\Sing(J,b)$ be the set of points $\xi\in W$ such that the
order of $J_{\xi}$ is $\geq b$. This is a closed subset of $W$. We
perform transformations of the following type: blow-up $W$ along a
suitable regular center $C\subset\Sing(J,b)$( which we will call a
{\it permissible center}), to get $W_{1}$. Let $E_{1}$ be
as in 1.1 and $J_{1}$ a certain transform of $J$ (see [11],
4.6, or
5.4 of this article for the precise definition of $J_{1}$). Then we
consider
$B_{1}=(W_{1},(J_{1},b),E_{1})$. This new basic object is called a
{\it permissible transform} of $B$.
\medskip

{\bf (5.2)} More explicitly, the idea is to define an
upper-semi-continuous
function $g$ from $W$ to a totally ordered set $\Lambda^{(d)}$ taking
finitely many values (the set $\Lambda^{(d)}$ depends only on the
dimension $d$ of $W$), so that
$C=\MaxBar{g}=\{\xi\in W\mid g(\xi)=\Max{g}\}$
is a regular subscheme of $\Sing(J,b)$. Take the
transform $B_{1}$ of $B$ with center $C$, define a function
$g_{1}:W_{1}\longrightarrow\Lambda^{(d)}$ analogous
to $g$ and so on.
All there should be naturally defined, and in this way we obtain a
sequence
of transformations of basic objects, say

$${\rm(5.2.1)} \quad
(W_{0},(J_{0},b),E_{0}) \longleftarrow\cdots\longleftarrow
(W_{r},(J_{r},b),E_{r})$$

The point is that, for any $ (W_{0},(J_{0},b),E_{0}) $, there
is an index $ N $, depending on this basic object, so that

$${\rm(5.2.2)} \quad
(W_{0},(J_{0},b),E_{0}) \longleftarrow\cdots\longleftarrow
(W_{N},(J_{N},b),E_{N})$$
is such that  the basic object $B_{N}=(W_{N},(J_{N},b),E_{N})$ satisfies
$\Sing(J_{N},b)=\emptyset$, i.~e. the order of $J_{N}$ is $<b$
at any point of $W_{N}$. The sequence (5.2.2) is called the { \it
resolution} of
  $ (W_{0},(J_{0},b),E_{0}) $ defined by the algorithm. The
functions $g_0, g_1, \ldots $ are called the {\it resolution
functions} that the algorithm attaches to the basic object $B_0$.
\medskip

{\bf (5.3) Remark.} Note that the resolution (5.2.2) is defined in
terms of
the functions $ g_{i}: W_i \longrightarrow \Lambda^{(d)}$. Setting
$C_i=\MaxBar{g_{i}}$, $C_i$ is a  subscheme of
$\Sing(J_i,b)$. Then, by taking restrictions,  we may think of $
g_{i} $ as a function defined on
$\Sing(J_i,b) \subset W_i $ (as is done in [11]). If one
follows this approach, the procedure comes
to an end when we reach the resolution (5.2.2) since
$\Sing(J_{N},b)=\emptyset$.
We shall later indicate how  $\Lambda^{(d)}$ and the functions
$ g_{i}: W_i \longrightarrow \Lambda^{(d)}$ are defined. We also refer
to [11] for properties of the resolutions (5.2.2) arising in
this way. Here we list some particularly important ones.
\smallskip

{\bf p0.} {\it Compatibility with open restrictions}. If
$ W_0 $ is replaced by an open subset, and all the other terms of
the basic
object are defined by restriction, then the resolution of the new
basic object and the corresponding resolution functions
$g_{i} $ are also defined by the natural restriction of (5.2.2) to such
open subset (after neglecting the arrows inducing isomorphism as in
1.2(vii)).
\smallskip

{\bf p1.} If $\xi\in\Sing(J_{i},b)$,
$i=0,\ldots,r-1$, and if $\xi\not\in C_{i}$ then
$g_{i}(\xi)=g_{i+1}(\xi')$ via the natural
identification of the point $\xi$ with a point $\xi'$ of
$\Sing(J_{i+1},b)$.
\smallskip

{\bf p2.} The resolution is obtained by transformations
with permissible centers $\MaxBar{g_{i}}$, for $i=0,\ldots,r-1$, and

$$\Max{g_{0}}>\Max{g_{1}}>\cdots>\Max{g_{r-1}}$$

{\bf p3.} If $J_{0}$ is the ideal of a regular pure
dimensional subvariety $X_{0}$, $E_{0}=\emptyset$
and $b=1$, then the function $g_{0}$ is constant.
\smallskip

{\bf p4.} For any $i=0,\ldots,r-1$, the closed set
$\MaxBar{g_{i}}$ is smooth, equidimensional and its
dimension is determined by the value $\Max{g_{i}}$.

The two problems of desingularization in (1.2) and of
resolution of basic objects are closely related, this will be
discussed in
(5.5).
\medskip

{\bf (5.4)} Here we shall be more explicit about some notions
introduced  in 5.1.

Let $B=(W,(J,b),E)$ be a basic object. The set $\Sing(J,b)$ of 5.1
can be described
more algebraically. Namely, it is the zero set of the $W$-ideal
denoted by  $\Delta_{W}^{b-1}(J)$, which is
the sheaf whose stalk at $\xi\in W$ is the ideal of
$\O_{W,\xi}$ generated by elements of $J_{\xi}$ as well as their
derivatives of order $<b$ (see 1.5 in [11]).
\smallskip

If $C\subset\Sing(J,b)$ is a permissible center (see 1.1.c) , in
5.1 we said that
  the transformation of $B$ with center $C$ is a certain basic
object $(W_{1}(J_{1},b),E_{1})$. The set $E_1$ is obtained as in
1.1.d, let us explain better how $J_{1}$
is obtained. If the center $C$ is irreducible,
the total transform $J\O_{W_{1}}$ can be written uniquely as
$$J\O_{W_{1}}=I(H)^{\nu}\bar{J}_{1}=I(H)^{b}J_{1}$$
for suitable ideals $\bar{J}_{1}$ and $J_{1}$ of $\O_{W_{1}}$, where
$I(H)$ is the ideal sheaf of the exceptional divisor $H$ and
$\nu=\nu_{\eta}(J)$, $\eta$ the generic point of $C$.
\smallskip

In this way we define $J_{1} \subset \O_{W_{1}}$. There is an
important relationship between $J_{1}$ and $\bar{J}_{1}$, namely:
$$ J_{1}=I(H)^{\nu-b}\bar{J}_{1}$$
and hence $ J_{1}=\bar{J}_{1}$ if $\nu=b $.
\smallskip

This construction can be extended to the case where $C$ is not
necessarily irreducible. One gets similar formal expressions, but now
the exponents must be interpreted as locally constant functions taking
integral values.
We write $B\longleftarrow B_{1}$ to indicate this transformation,
if the center $C$ is clear from the context.
\smallskip

If we iterate this process, starting from
$B=B_{0}=(W_{0},(J_{0},b),E_{0})$ we
get a sequence of transformations of basic objects
$${\rm(5.4.1)} \quad
(W_{0},(J_{0},b),E_{0}) \longleftarrow\cdots\longleftarrow
(W_{r},(J_{r},b),E_{r})$$
and for each $i=0,\ldots,r$ an expression
$${\rm(5.5.2)} \quad
J_{i}=I(H_{r+1})^{a_{1}}\cdots I(H_{r+i})^{a_{r}}\bar{J}_{i}$$
where $H_{j}$ is the strict transform of the exceptional divisor of
$W_{j}\longrightarrow W_{j-1}$ (see 4.8 in [11] for more details).
\smallskip

We let $\nu_{i}:W_{i}\longrightarrow{\Bbb Z}$ denote the order
function corresponding of the ideal $\bar{J}_{i}\subset\O_{W_{i}}$
(the biggest power of the maximal ideal of $ {\cal O}_{W_{i},x} $
  containing  $ (\bar{J}_{i})_{x} $).
\medskip

{\bf (5.5)} Let us  indicate how the existence of an algorithm of
resolution of basic objects implies one of strong principalization
(1.2).
Assume an id-triple $(W,{\cal I}, E)$ (as in 1.1.1) is given. To
obtain a principalization thereof,
consider
$B_0= (W_{0},(J_{0},1),E_{0})$,
where $W_{0}=W$, $J_0={\cal I}$, $E_0=E$. Let (5.2.2) be the resolution
of $B_0$ defined by the algorithm (applied to $B_0$) and note that
\smallskip

{\bf (i)} for each index $ r $
$$ J_{0}{\cal O}_{W_r}=I(H_{r+1})^{c_{1}}\cdots I(H_{r+i})^{c_{r}}
J_{r} $$
for suitable integers $ c_i $.
\smallskip

{\bf (ii)} $ J_{0}{\cal O}_{W_N}=I(H_{r+1})^{c_{1}}
\cdots I(H_{r+N})^{c_{r}} {\cal O}_{W_N} $
\smallskip

Now we can easily check that all conditions in 1.2 hold, by setting
$I^{(d)}=\Lambda^{(d)}$ and $h_i=g_i$, for all $i$, where  $g_0,
g_1, \ldots $ are the resolution functions that the algorithm
attaches to the basic object $B_0$.
\medskip

{\bf (5.6)}
Recall that a sequence (5.4.1) is defined by centers
$C_i\subset\Sing(J_i,b)$ (5.1). Assume inductively that
$C_{i}\subset\Max\nu_{i}$ (here we mean $\nu_{i}$ restricted to the
closed set $\Sing(J_i,b)$) $i=0,1,\ldots,r-1$. In this setting we attach
to each term
$B_{i}=(W_{i},(J_{i},b),E_{i})$
a function $t_{i}:\Sing(J_i,b) \longrightarrow{\Bbb
Q}\times{\Bbb Z}$ as
follows.
Let $\word_{i}(\xi)=\nu_{i}(\xi)/b$. Under our assumption on the
centers $C_{i}$, one can show (see 4.12.1 in [11]) that
$\Max\word_{i}\geq\Max\word_{i+1}$
for all $i$.
For a given index $i=0,1,\ldots,r$, assume that $i_{0}$ is such that
$$\Max\word_{i_{0}-1}>\Max\word_{i_{0}}=
\Max\word_{i_{0}+1}=\cdots=\Max\word_{i}$$
Let $E_{i}^{-}$ consist of the hypersurfaces of $E_{i}$ which are
strict transforms of the ones of $E_{i_{0}}$, and $E_{i}^{+}$ the
remaining ones. Set for $\xi\in W_{i}$
$$n_{i}(\xi)=\#\{H\in E^{-}_{i}\mid \xi\in H\}$$
We set $t_{i}(\xi)=(\word_{i}(\xi),n_{i}(\xi))$. Then $t_{i}$ is
an upper-semi-continuous function, when ${\Bbb Q}\times{\Bbb Z}$ is
ordered lexicographically (see 4.15 in [11]).
The sequence (5.4.1) will be called $t$-{\it permissible} if
$C_{i}\subset\Max{t_{i}}$ for $i=0,1,\ldots,r-1$.
In this case, we have that (see 4.15 in [11]):
$$\Max{t_{i}}\geq\Max{t_{i+1}} \qquad \forall\ i$$
The key step in the resolution of basic objects is to drop the
maximum of the function
$\word$. In fact, if we achieve a situation where $\Max\word_{N}=0$
for some index $N$, then the
resolution of the basic object $B_{N}$ becomes a rather simple
combinatorial problem
(section~5 in [11]). Note that in this case $\bar{J}_{N}=\O_{W_{N}}$.
\medskip

{\bf (5.7)} We recall some conditions, important in  inductive
arguments to be used in the next chapter.
Let $(W,(J,b),E)$ be a basic object.
Fix a point $\xi\in\Sing(J,b)$.
We say that the basic object satisfies conditions (LC) and (IA) at
$\xi$ if there exists a smooth hypersurface
$Z$ of $W$  defined on a neighborhood of $\xi$ such that:
\smallskip

(LC) {\bf Local condition.} $I(Z)\subset\Delta^{b-1}(J)$.
\smallskip

(IA) {\bf Inductive assumption.} $Z$ is transversal to all the
hypersurfaces of $E$.
\smallskip
We remark that conditions (LC) and (IA) are stable under
permissible transformations. This means the following. If
$(W_{1},(J_{1},b),E_{1})\longrightarrow(W,(J,b),E)$ is a
transformation of basic objects and $Z_{1}\subset W_{1}$ is the
strict transform of $Z$, then (LC) and (IA) also hold for $Z_{1}$ and
$(W_{1},(J_{1},b),E_{1})$.
The result about (LC) follows from ``Giraud's Lemma''
(see lemma~6.6 in [11])
  and that for for (IA)  is well known .
\smallskip
The resolution of basic objects will be reduced, using the function
$t$ (see 5.6), to the special case
where condition (LC) and (IA) hold. But this last case can be
understood as that of a basic object in a space whose dimension has
dropped by one. Note also that if condition (LC) holds, then at any
point
$\xi\in\Sing(J,b)$ the order of $J_{\xi}$ is $b$. The following
result indicates how this is done.

In it, the following notation is used. If $Y$ is a subscheme of
$W$, $R(1):=R(1)(Y)$ denotes the
union of the irreducible components of $Z$ of codimension one;
$C(J)$ is a certain sheaf of ideals on $Z$ (as above), called the
{\it coefficient ideal} of $J$. See 9.3 or exercise 13.1 of
[11] for its definition and basic properties. But if $b=1$,
$C(J)=J{\cal O}_Z$ (for instance, if $J$ is the ideal of a regular
scheme $X \subset W$.)
\medskip

{\bf (5.8)} {\bf Inductive Lemma}
([11], Lemma~6.12).
Let $B=(W,(J,b),E)$ be a basic object such that (LC) and (IA) hold
at any point of $\Sing(J,b)$. Then:
\smallskip
{\bf (a)} If $R(1)(\Sing(J,b))\neq\emptyset$ then
then $R(1)$ is regular, open and closed in $\Sing(J,b)$ and has normal
crossings with $E$.
Consider the transformation
$(W,(J,b),E)\longleftarrow (W_{1},(J_{1},b),E_{1})$
with center $C=R(1)$, then $R(1)(\Sing(J_{1},b))=\emptyset$.
\smallskip
{\bf (b)} If $R(1)(\Sing(J,b))=\emptyset$, we define a basic object
$B_{Z}=(Z,(C(J),b!),E_{Z})$, where $C(J)$ is
the coefficient ideal of $J$
and $E_{Z}$ consists of the intersection of
hypersurfaces in $E$ with $Z$.
This basic object $B_{Z}$ satisfies the following properties:
\smallskip
{\bf (b1)} $\Sing(J,b)=\Sing(C(J),b!)$.
\smallskip
{\bf (b2)} For any sequence of transformations, say
$$(W,(J,b),E)\longleftarrow\cdots\longleftarrow
    B_{s}=(W_{s},(J_{s},b),E_{s})$$
with centers $C_{i}\subset\Sing(J_{i},b)$ then
there is a sequence of transformations
$$(Z,(C(J),b!),E_{Z})\longleftarrow\cdots\longleftarrow
(Z_{s},(C(J)_{s},b!),(E_{Z})_{s})$$
such that
$\Sing(C(J)_{i},b!)=\Sing(J_{i},b)$ and with the same centers $C_{i}$.
\medskip

{\bf (5.9)} To resolve, inductively, basic objects, in [11]
one uses a certain (locally defined) auxiliary basic object, denoted
$B''$ (or $B_r''$, if appropriate)), satisfying conditions (IC) and
(IA). But it is necessary
to construct first another basic object, denoted $B'$, related to
the function w-ord (the singular locus of $B'$ is the
closed set $\MaxBar\word$ and the resolution of $B'$ implies that the
value $\Max\word$ has dropped). Unfortunately, in general $B'$ only
satisfies condition (LC)
and not condition (IA). This $B'$ is defined as follows.
\smallskip

Consider a sequence of transformations of basic objects,
$$(5.9.1) \qquad (W_{0},(J_{0},b),E_{0})
\longleftarrow\cdots\longleftarrow (W_{r},(J_{r},b),E_{r})$$
such that the centers $C_{i}\subset\Max\word_{i}$.
Set $b_{r}=b\cdot\Max\word_{r}$ and assume that $\Max\word_{r}$ is
strictly positive. Then,
there exists a basic object $B'_{r}=(W_{r},(J'_{r},b'),E_{r})$ such that
\smallskip

{\bf (1)} Condition (LC) holds locally for $(W_{r},(J'_{r},b'),E_{r})$.
\smallskip

{\bf (2)} $\MaxBar\word_{r}=\Sing(J'_{r},b')$ and this condition is
stable after transformation, see below.
\smallskip

{\bf (3)} Any sequence of transformations
$$(W_{r},(J'_{r},b'),E_{r})\longleftarrow\cdots\longleftarrow
(W_{s},(J'_{s},b'),E_{s})$$
induces an extension of the sequence (5.9.1) with same centers
$C_{i}\subset\Sing(J'_{i},b')$,
$$(W_{r},(J_{r},b),E_{r})\longleftarrow\cdots\longleftarrow
(W_{s},(J_{s},b),E_{s})$$
with
\smallskip

{\bf (3a)} $\MaxBar\word_{i}=\Sing(J'_{i},b')$ for $r\leq i<s$ and
$$\Max\word_{r}=\cdots=\Max\word_{s-1}$$
{\bf (3b)} $\Sing(J'_{s},b')=\emptyset$ if and only
$\Max\word_{r}>\Max\word_{s}$; and if $\Max\word_{r}=\break\Max\word_{s}$
then $\MaxBar\word_{s}=\Sing(J'_{s},b')$.
\smallskip

The definition of the basic object $B_r ' =
(W_{r},(J'_{r},b'),E_{r})$ is the
following (see [11], 9.5.4):
$$J'_{r} =\left\{
{\bar{J}_{r}
\atop
\bar{J}_{r}^{b-b_{r}}+
   \left(I(H_{1})^{a_{1}}\cdots I(H_{r})^{a_{r}}\right)^{b_{r}}}
\hskip 1cm
{{\rm if} \atop {\rm if}}
\hskip 1cm
{{b_{r}\geq b} \atop {b_{r}<b}}
\right.$$
$$b'=\left\{
{b_{r} \atop b_{r}(b-b_{r})}
\hskip 1cm
{{\rm if} \atop {\rm if}}
\hskip 1cm
{{b_{r}\geq b} \atop {b_{r}<b}}
\right.$$
\medskip

{\bf (5.10)} Now we are in position to construct the basic object
$B''$, which describes the set
$\MaxBar{t}$ and  satisfies conditions (LC) and (IA). (Hence
we may apply the inductive lemma~5.8 to this basic object). This is the
procedure to achieve resolution of basic objects: we use the function
$t$ to obtain a basic object satisfying conditions (LC) and (IA) and
then we use induction on the dimension.
\smallskip

Consider a $t$-permissible sequence of transformations of basic objects
as in (5.4.1),
$$(W_{0},(J_{0},b),E_{0})\longleftarrow\cdots\longleftarrow
(W_{r},(J_{r},b),E_{r})$$
There exists a basic object $B''_{r}=(W_{r},(J''_{r},b''),E''_{r})$
such that
\smallskip

{\bf (1)} Conditions (LC) and (IA) hold for
$(W_{r},(J''_{r},b''),E''_{r})$, locally at any point
$\xi\in\Max{t_{r}}$.
\smallskip

{\bf (2)} $\MaxBar{t_{r}}=\Sing(J''_{r},b'')$.
\smallskip

{\bf (3)} Any sequence of transformations
$$(W_{r},(J''_{r},b''),E''_{r})\longleftarrow\cdots\longleftarrow
(W_{s},(J''_{s},b''),E''_{s})$$
induces an extension of the sequence (5.9.1) with same centers
$C_{i}\subset\Sing(J''_{i},b'')$,
$$(W_{r},(J_{r},b),E_{r})\longleftarrow\cdots\longleftarrow
(W_{s},(J_{s},b),E_{s})$$
with
\smallskip

{\bf (3a)} $\MaxBar{t_{i}}=\Sing(J''_{i},b'')$ for $r\leq i<s$ and
$$\Max{t_{r}}=\cdots=\Max{t_{s-1}}$$
\smallskip

{\bf (3b)} $\Sing(J''_{s},b'')=\emptyset$ if and only
$\Max{t_{r}}>\Max{t_{s}}$; and if $\Max{t_{r}}=\Max{t_{s}}$
then $\MaxBar{t_{s}}=\Sing(J''_{s},b'')$.
\medskip

The definition of the basic object $(W_{r},(J''_{r},b''),E''_{r})$
is, locally at any point $\xi$ (see 9.5.7 in [11]), as follows:
$$J''_{r}=(J'_{r})+
I(H_{i_{1}})^{b''}+\cdots+I(H_{i_{N}})^{b''}$$
where $H_{i_{1}},\ldots,H_{i_{N}}$ are the hypersurfaces of $E_{r}^{-}$
passing through the point $\xi$,
$b''=b'$ and $E''_{r}=E^{+}_{r}$.
\medskip

{\bf (5.11)}
Finally, we are in condition to describe more precisely the functions
$g_{i}$ of 5.2. We define them by induction on the dimension $d$ of the
ambient space $W$. In all cases, the values will be in a totally
ordered set $\Lambda^{(d)}$ with a first element $0_{d}$. We will
have $g_{i}(\xi)=0_{d}$ if $\xi\not\in\Sing(J_{i},b)$. So we may
restrict ourselves to points $\xi\in\Sing(J_{i},b)$.
For $d=1$ we set $\Lambda^{(1)}=\{0_{1}\}\cup({\Bbb Q}\times{\Bbb 
Z})\cup\infty_{1}$, where $\Lambda^{(1)}$ is ordered 
lexicographically, $0_{1}$ is the first element and $\infty_{1}$ is 
the last one,
and set
$$g_{i}(\xi)=
\left\{
{t_{i}(\xi)\atop \Gamma_{i}(\xi)}
\hskip 1cm
{{\rm if} \atop {\rm if}}
\hskip 1cm
{\word_{i}(\xi)>0 \atop \word_{i}(\xi)=0}
\hskip 1cm
{\hbox{(see (5.6))}
\atop
\hbox{(see [11] section 5})}
\right.$$ 
\medskip

Assume $d>1$ and that the resolution functions are defined when
dimension is $<d$.
Set
$$\Lambda^{(d)}=
\left(\left(
	\left({\Bbb Q}\times{\Bbb Z}\right)\cup
	\left({\Bbb Z}\times{\Bbb Q}\times{\Bbb Z}^{\Bbb N}\right)
\right)
\times\Lambda^{(d-1)}\right)\cup\{0_{d},\infty_{d}\}$$
naturally ordered as in 9.6 of [11]; $0_{d}$ is the first element
and $\infty_{d}$ is the last one.

For dimension $d$ the functions $g_{i}$ take values in $\Lambda^{(d)}$
and are defined inductively on the length $r$ of the sequence
involved, as follows.

Fix $r=0$, i.e. the sequence consists in just one term
$B=(W,(J,b),E)$. Let $\xi \in \Sing(J,b)$ and we may consider a
suitable neighborhood of $\xi\in W$ such that $\xi\in\MaxBar{t_{0}}$.
In principle there are three possible cases, for the point $\xi$:
\smallskip

{\bf (1)} $\word_r(\xi)=0$.
\smallskip

{\bf (2)} $\xi\in R(1)(\MaxBar{t_r})$
(the union of one-codimensional components of $\MaxBar{t_r}$) and 
$\word_{r}(\xi)>0$.
\smallskip

{\bf (3)} $\xi\not\in R(1)(\MaxBar{t_r})$ and $\word_r(\xi)>0$.
\smallskip

But, in fact, case (1) cannot occur when $r=0$, so we consider the
other two.
\smallskip

In case (2) we set $g_{r}(\xi)=(t_{r}(\xi),\infty)$.
\smallskip

In case (3), consider the associated basic object $B''$ of 5.10.
Then by 5.10(1) near $\xi$ there is a smooth hypersurface $Z_{0}$
satisfying conditions (LC) and (IA), and a basic object
$B_{Z_{0}}=(Z_{0},(C(J_{0}),b''!),E_{Z_{0}})$.
Since $\dim{Z_{0}}=d-1$, by induction there
is an associated resolution function
$g'_{0}:Z_{0}\longrightarrow\Lambda^{(d-1)}$. Set
$g_{0}(\xi)=(t_{0}(\xi),g'_{0}(\xi))$.
It can be proved that this value is independent
of the choice of $Z_{0}$ (see [11], Section 7).
\smallskip

Now, assuming the definition given for sequences of length $<r$,
consider one of length $r$, say 5.4.1. By induction, it
suffices to define $g_{r}:W_{r}\longrightarrow\Lambda^{(d)}$.
Let $\xi_{r}\in\Sing(J_{r},b)$ and consider a
suitable neighborhood of $\xi_{r}\in W_{r}$ such that
$\xi_{r}\in\MaxBar{t_{r}}$.
Again there are three cases, as for $r=0$, for the point $\xi_{r}$.
In case (1), necessarily
$J_{r}=I(H_1)^{a_1}\ldots I(H_m)^{a_m}$ (where $E_{r}=\{H_1,\ldots,H_m
\}$). We write
$g_r(\xi_{r})=(\Gamma_r(\xi_{r}),\infty)$ (see [11], section~5);
$\Gamma _r $ is a certain function from $\Sing(J_{r},b)$ to
${\Bbb Z}\times{\Bbb Q}\times{\Bbb Z}^{\Bbb N}$.
The definition in case (2) is
exactly like that for $r=0$.
In case (3), set $\xi_{i}\in W_{i}$ the image of $\xi_{r}$ via the  
morphism
$W_{r}\longrightarrow W_{i}$.
Let $r_{0}$ be the smallest index so that
$t_{r}(\xi_{r})=t_{r_{0}}(\xi_{r_{0}})$.
In a suitable neighborhood of
$\xi_{r_{0}}$ we may assume that
$\Max{t_{r_{0}}}=t_{r_{0}}(\xi_{r_{0}})$ and we may consider the
associated basic objects
$B''_{r_{0}}, B''_{r_{0}+1},\ldots, B''_{r}$ (5.10).
There are sequences of transformations of basic
objects
$$B''_{r_{0}}\longleftarrow B''_{r_{0}+1}
\longleftarrow\cdots\longleftarrow B''_{r}$$
and
$$B_{Z_{r_{0}}}\longleftarrow B_{Z_{r_{0}+1}}
\longleftarrow\cdots\longleftarrow B_{Z_{r}}$$
where each $Z_{i}$ is a hypersurface defined near $\xi_{i}$ in
$W_{i}$, satisfying (IA) and (LC). By induction on the dimension
$d$ there are resolution
functions for the last sequence, say
$g'_{i}:Z_{i}\longrightarrow\Lambda^{(d-1)}$, for
$i=r_{0},r_{0}+1,\ldots,r$ and we set
$$g_{r}(\xi)=(t_{r}(\xi),g'_{r}(\xi))$$
This is a well defined function (see [11], Section~7.)
\bigskip

{\bf 6. PROOF OF CONDITION 1.13 (5) IN THEOREM 1.16.}
\medskip

In this chapter we shall see that the algorithm of 1.16 (that is,
the one discussed in chapter 5, referred henceforth as {\it the
algorithm}) satisfies condition 1.13.(5), i.e., it is a good
principalization algorithm. It will be easily seen (cf. 6.3) that
this result rapidly follows from a similar one about basic objects
(Theorem 6.4). To make this precise, first we need a notion of
compatibility, analagous to that of 1.11.
\medskip

{\bf  6.1.}{\it Compatibility.}
Let $B_{0}=(W_{0},(J_{0},b),E_{0})$ be a basic object and let
$W^{(t)}_{0}\subset W_{0}$ be a pure dimensional smooth closed
subscheme. Let 5.2.2 be the resolution sequence that the algorithm
attaches to $B_0$.

We shall define the notion: ``$W^{(t)}_{0}$ is r-compatible with
$B_0$, where $0 \leq r \leq N$.
\smallskip

We say that $W^{(t)}_{0}$ is $0$-compatible with $B_{0}$
if for any hypersurface $H \in E_{0}$,  $H$ and $W^{(t)}_{0}$ meet
transversally and $H\cap W^{(t)}_{0}$ is a
hypersurface in $W^{(t)}_{0}$; furthermore we require that
$E^{(t)}_{0}=\{H\cap W^{(t)}_{0}\mid H\in E_{0}\}$ is a set of smooth
hypersurfaces in $W^{(t)}_{0} $ having only normal crossings.
Note that in such case
$B^{(t)}_{0}=(W^{(t)},(J^{(t)}_{0},b),E^{(t)}_{0})$
is also a basic object, where $J^{(t)}_{0}=J_{0}\O_{W^{(t)}_{0}}$
is a non
zero sheaf of ideals.
\smallskip

Thus, we have two basic objects, namely $ B_{0}$ and  $ B^{(t)}_{0}$.
The algorithm will attach to each a resolution. If $r > 0$, consider the
  first $r$ steps of each resolution sequence, say

$${\rm (6.1.1)}\qquad
B_{0}\longleftarrow B_{1}\longleftarrow \cdots \longleftarrow B_{r}$$
and
$${\rm (6.1.2)}\qquad
B^{(t)}_{0}\longleftarrow B^{(t)}_{1}\longleftarrow \cdots
\longleftarrow B^{(t)}_{r}$$
respectively.
Write
$B_{i}=(W_{i},(J_{i},b),E_{i})$ and
$B^{(t)}_{i}=(W^{(t)}_{i},(J^{(t)}_{i},b),E^{(t)}_{i})$, let
$C_{i}=\MaxBar{g_{i}}$ and $C^{(t)}_{i}=\MaxBar{g^{(t)}_{i}}$ be the
$i$-centers of $B_{i}$ and $B^{(t)}_{i}$, defined by the resolution
functions, respectively,
$i=0,\ldots,r-1$.
Then we say that  $W^{(t)}_{0}$ is r-compatible with $ B_0 $ if
the following two conditions hold:
\smallskip

(C1) each  $C_{i}$ is transversal to $W^{(t)}_{i}$ and,

(C2) $C^{(t)}_{i}=C_{i}\cap W^{(t)}_{0}$,
$i=0,\ldots,r-1$.
\smallskip
Here we make use of the fact that if $W_{i+1}\longrightarrow W_{i}$
is the
blowing-up
with center $C_{i}$ and $W^{(t)}_{i+1}\longrightarrow W^{(t)}_{i}$
the blowing-up
with center $C^{(t)}_{i}$, then $W^{(t)}_{i+1}$ is naturally
identified with the
strict transform of $W^{(t)}_{i}$ in $W_{i+1}$. In particular it
makes sense to
require the transversality of $C_{i+1}$ and $W^{(t)}_{i+1}$.

If r-compatibility holds, then (6.1.1) and (6.1.2) are
strongly related. In fact, one can check that
$J^{(t)}_{i+1}=J_{i+1}\O_{W^{(t)}_{i+1}}$ and
$E^{(t)}_{i+1}=\{H\cap W^{(t)}_{i+1}\mid H\in E_{i+1}\}$ for each
index $ i
\leq r $.
\medskip

{\bf (6.2) Remark.} We may naturally define the notion of
{\it family of basic objects} as a   system $(W,(J,b),E,T,\Pi)$ where
$(W,(J,b),E)$ is a basic object $\Pi:W\longrightarrow T$ is a smooth
morphism of regular schemes, so that for each $t\in T$ the fiber
$(W^{(t)},(J^{(t)},b),E^{(t)})$ (defined by taking
$W^{(t)}=\Pi^{-1}(t)$, $J^{(t)}=J\O_{W^{(t)}}$ and $E^{(t)}$ the
  intersection of $E$ with $W^{(t)}$) is a basic object, for all $t$
in $T$.
In particular, we are requiring that $E^{(t)}$ should also
have normal crossings, which is a strong condition. Then, the most
natural example of the situation described in 6.1 would be the case
where $W_0^{(t)}$ is the space of a closed fiber of a family of
basic objects. As in 1.11, we could define the notion ``the family
is r-compatible with the algorithm''. This won't be relevant in our
present work, so we leave the details to the interested reader.
\medskip

{\bf (6.3) Remark.} Let
${\cal T}_0=(W_0,J_0,E_0)$ be an id-triple in $S_0$ (cf 1.1 and
1.15) As explained in 5.5, the p-sequence of ${\cal T}_0$ is
obtained from the resolution sequence 5.2.2 of the basic object
  $B_0=(W_0,(J_0,1),E_0)$. Assume $W^{(t)}$ is a reduced  closed
subscheme of  $W_0$ which is r-compatible with ${\cal T}_0$ (in the
sense of 1.10, relative to the algorithm). Then it is
straightforward to see that
$W^{(t)}$ is r-compatible with $B_0$ in the sense of 6.1. Since the
resolution functions that appear in the p-sequence of ${\cal T}_0$
and  those of the the resolution sequence 5.4.1 of $B_0$ are the
same, from the discussion in 1.11 it follows that to show 1.13 (5)
it suffices to show the following result:
\proclaim
(6.4)  Theorem.
Let $B_{0}=(W_{0},(J_{0},b),E_{0})$ be a basic object,
$W_{0}^{(t)}\subset W_{0}$ smooth closed subscheme and $r\geq 0$.
Assume that $W^{(t)}_{0}$ is
$r$-compatible with $B_{0}$. If $\xi_{r}\in W^{(t)}_{r}$ then
\smallskip
(a) $g_{r}(\xi_{r})\leq g_{r}^{(t)}(\xi_{r})$.
\smallskip
(b) $g_{r}(\xi_{r})=g_{r}^{(t)}(\xi_{r})$ if and only if
$S_{r,\alpha_{r}}$ is transversal to $W_{r}^{(t)}$ at $\xi_{r}$,
where $\xi_{r}\in S_{r,\alpha_{r}}$ and $\alpha_{r}=g_{r}(\xi_{r})$
(the notation is similar to that of 1.13 and 1.14.)

{\bf (6.5)} This result will be a consequence of Propositions 6.14
through 6.17. The content of these propositions, taken together, is
the same  as that of 6.4, but the hypotheses are stronger. Namely,
certain ''extra assumptions" are made. However, in 6.18 we'll see
that that, except for some rather trivial special cases, we may
assume the validity of these extra sssumptions, obtaining in this
way a proof of 6.4. Before discussing the cited propositions, let us
discuss the ``extra assumptions''.
\medskip

{\bf (6.6)} Our assumptions and notation are those of Theorem 6.4.
  Consider the
sequence (6.1.1), fix $\xi_{r}\in W^{(t)}_{r}$ and let $\xi_{i}\in W_{i}$
denote the image of $\xi_{r}$
in $W_{i}$. Thus, we have a sequence of points
$\xi_{0},\xi_{1},\ldots,\xi_{r}$

The {\it local extra assumptions} are the following:

(1) $g_{i}(\xi_{i})=g^{(t)}_{i}(\xi_{i})$ (say $=\alpha_{i}$), for
$i=0,\ldots,r-1$.

(2) The stratum $S_{i,\alpha_{i}}$ is transversal to $W^{(t)}_{i}$ at $
\xi_{i} $, and $S_{i,\alpha_{i}} \cap W^{(t)}_{i}=S^{(t)}_{i,\alpha_{i}}$
locally at $ \xi_{i} $, for $i=0,\ldots,r-1$, where $
S^{(t)}_{i,\alpha_{i}}$ denotes the stratum of $ B^{(t)}_i $ containing $
\xi_{i} $ in (6.1.2).
\medskip
{\bf (6.7)} Recall that given a sequence of permissible
transformations of basic objects 5.4.1, we have an expression
involving $W_i$-ideals:
$$J_{i}=I(H_{1})^{a_{1}}\cdots I(H_{i})^{a_{i}}\bar{J}_{i}, ~
i=0,\ldots,r-1,r$$.

Under the assumtions of 6.6 (including (1) and (2) therein), we
claim that:
\smallskip
(3) locally at $ \xi_i $, and
for each index $i=0,\ldots,r-1,r$, the corresponding expression of
$J^{(t)}_{i}$ for
$i=0,\ldots,r$ will be
$$J^{(t)}_{i}= I(H^{(t)}_{1})^{a_{1}}\cdots
I(H^{(t)}_{i})^{a_{i}}\bar{J}^{(t)}_{i},$$
where $H^{(t)}_{j}=H_{j}\cap W^{(t)}_{i}$, $j=1,\ldots,i$ and
$\bar{J}^{(t)}_{i}=\bar{J}_{i}\O_{W^{(t)}_{i}}$.
\smallskip
In fact, note first that for $ r=0$ condition (3) holds since
$J_{0}=\bar{J}_{0}$ and $J^{(t)}_{0}=\bar{J}^{(t)}_{0}$.

So, proceed by induction and assume that condition (3) holds for
$i=0,\ldots,r-1$. If $\xi_{r-1}\not\in C_{r-1}$ the (3) holds since, in
this case, we may identify $\xi_{r-1}$ with $ \xi_{r}$. So assume that
$\xi_{r-1}\in C_{r-1}$ and recall from 5.6.2 that the expression of
$J_{r}$ (resp.  $J^{(t)}_{r}$) is defined in terms of the expression of
$J_{r-1}$ (resp.  $J^{(t)}_{r-1}$), the number of hypersurfaces
$H_{j}$ (resp.  $H^{(t)}_{j}$) included in $C_{r-1}$ (resp.
$C^{(t)}_{r-1}$) and the order of $\bar{J}_{r-1}$ (resp.
$\bar{J}^{(t)}_{r-1}$) along $C_{r-1}$ (resp. $C^{(t)}_{r-1}$). Let us
check this last point, and assume, in addition, that $\bar{J}_{r-1}$ is a
proper ideal ( that $\word_{r-1}(\xi_{r-1})>0$), otherwise (3) follows
easily by induction and (2). Since $\word_{r-1}(\xi_{r-1})$ is a
coordinate of
$ g_{r-1}(\xi_{r-1})$ it turns out that
$$ \nu_{\xi_{r-1}}(\bar{J}_{r-1})=\nu_{C_{r-1}}(\bar{J}_{r-1})$$

where $\nu_{C_{r-1}}(\bar{J}_{r-1})$ denotes the order of $J_{r-1}$
at the generic point of the component of $C_{r-1}$ which contains
$\xi_{r-1}$.

The same argument applies for $\bar{J}^{(t)}_{r-1}$ and
$C^{(t)}_{r-1}$, so
$$ \nu_{\xi_{r-1}}(\bar{J}^{(t)}_{r-1})=
\nu_{C^{(t)}_{r-1}}(\bar{J}^{(t)}_{r-1})$$

Finally, since $g_{r-1}(\xi_{r-1})=f^{(t)}_{r-1}(\xi_{r-1})$ (by (1)), it
follows that
$ \nu_{C_{r-1}}(\bar{J}_{r-1})=\nu_{C^{(t)}_{r-1}}(\bar{J}^{(t)}_{r-1})$.
This proves the validity of (3), when (1) and (2) hold.
\smallskip

Before we  discuss the important propositions announced in 6.5, we
need some preliminary results.
\proclaim
(6.8) Lemma.
Let $W$ be a smooth scheme and let $J\subset\O_{W}$
be an ideal. Then:
\smallskip
(a) At any point $\xi\in W^{(t)}$, the order of $J$ at $\xi\in W $
is smaller
or equal than the order of $J^{(t)}$ at the same point $\xi\in W^{(t)}$ .
\smallskip
(b) Set $b$ the maximum value of the order of $J$.
Consider now two smooth  subschemes of $W$, say $C$ and $W^{(t)}$,
  such that
$C\subset\Sing(J,b)$ (so the order of $J$ along $C$ is
constant and equal to $b$), and set $J^{(t)}=J\O_{W^{(t)}}$.
If $C$ is transversal to $W^{(t)}$ then for any
$\xi\in C\cap W^{(t)}$ the order of $J^{(t)}$ is also $b$.

{\it Proof:}
Denote by ${\cal M}$ the Maximal ideal of $\O_{W,\xi}$, $\bar{{\cal M}}$
the Maximal ideal of $\O_{W^{(t)},\xi}$.
The order of $J$ at the point $\xi$ is the biggest integer
$d$ such that $J_{\xi}\subset{\cal M}^{d}$.
It is clear that $J^{(t)}_{\xi}\subset\bar{\cal M}^{d}$, so that the
order of $J^{(t)}$ at $\xi$ is greater or equal than the order of $J$
and (a) is proved.

For (b), let $\xi\in C\cap W^{(t)}$. Transversality ensures that
there exists a regular system of
parameters $x_{1},\ldots,x_{n}$ of $\O_{W,\xi}$ such that
$I(C)_{\xi}={\cal P}=(x_{1},\ldots,x_{r})$ and
$I(W^{(t)})_{\xi}=(x_{s},\ldots,x_{n})$.
Since $C\subset\Sing(J,b)$, the order of the ideal $J_{\xi}$ at the
Maximal ideal ${\cal M}$ is also the biggest integer $b$ such that
$J_{\xi}\subset{\cal P}^{b}$.
So there is an element $f\in J_{\xi}$ such that
$f\in J_{\xi}\setminus{\cal M}^{b+1}$.
As $J_{\xi}\subset{\cal P}^{b}$ we may find an expression
$$f=\sum_{i_{1}+\cdot+i_{r}=b}a_{i_{1},\ldots,i_{r}}x_{1}^{i_{1}}\cdots
x_{r}^{i_{r}}$$
where $a_{i_{1},\ldots,i_{r}}\in\O_{W,\xi}$; and there must be indices
$i_{1},\ldots,i_{r}$ such that $a_{i_{1},\ldots,i_{r}}$ is a unit in
$\O_{W,\xi}$.

Consider the class $\bar{f}$ of $f$ in $\O_{W^{(t)},\xi}$ and note
that the above expression still holds considering the corresponding
classes. Furthermore, $\bar{a}_{i_{1},\ldots,i_{r}}$ is a unit in
$\O_{W^{(t)},\xi}$ for some index $i_{1},\ldots,i_{r}$.
We have proved that
$\bar{f}\in J^{(t)}_{\xi}\setminus{\cal M}^{b+1}$, and hence the order of
$J^{(t)}$ at $\xi$ is also $b$.

\proclaim
(6.9) Lemma.
Let $W$ be a smooth variety, $J\subset\O_W$ and
$W^{(t)}$ be a smooth subscheme of $W$. Set $J^{(t)}=J\O_{W^{(t)}}$.
Fix a point $\xi\in W^{(t)}$.
Assume that $\nu_J(\xi)=\nu_{J^{(t)}}(\xi)=b$.
\smallskip
Then there exists, locally at $\xi$, a
smooth hypersurface $Z$ of $W$ such that:
\smallskip
(1) $Z$ is transversal to $W^{(t)}$ at $\xi$.
\smallskip
(2) $I(Z)\subset\Delta_W^{b-1}(J)$ (locally at $\xi$).
\smallskip
(3) $I(Z^{(t)})\subset\Delta_{W^{(t)}}^{b-1}(J^{(t)})$ (locally at $\xi$)
where $ Z^{(t)}=Z \cap W^{(t)}$.

{\it Proof:}
One can check that
$\Delta_{W^{(t)}}(J^{(t)})\subset\Delta_W(J)\O_{W^{(t)}}$ locally
at $\xi$;
and hence, by
induction, that
$\Delta_{W^{(t)}}^{b-1}(J^{(t)})\subset\Delta_W^{b-1}(J)\O_{W^{(t)}}$.
Note that by assumption the order of the ideal
$\Delta_{W^{(t)}}^{b-1}(J^{(t)})$ at $\xi$ is one, so there exists
$\bar{x}_1\in\Delta_{W^{(t)}}^{b-1}(J^{(t)})_{\xi}$ of order one.
Let $x_1\in\Delta_W^{b-1}(J)_{\xi}$ be such the image in
$\O_{W^{(t)},\xi}$ is $\bar{x}_1$.
The hypersurface $Z$ defined by the equation $x_1$ (locally at $\xi$)
fulfills the required conditions.

\proclaim
(6.10) Lemma.
Fix $\xi_{r}\in W^{(t)}_{r} \subset W_{r} $ and assume that
conditions (1), (2) of  6.6 and (3) of 6.7 hold at
$\xi_{0},\xi_{1},\ldots,\xi_{r}$.
In addition assume that:
\smallskip
(i) $\word_r(\xi_{r})>0$ (so that
$g_r(\xi_{r})=(t_r(\xi_{r}),g'_r(\xi_{r}))$ (see 5.11 case (3))).
\smallskip
(ii) $t_r(\xi_{r})=t_r^{(t)}(\xi_{r})$.
\smallskip
After restriction to a suitable neighborhood of $\xi_{r}$ we may
assume that
$t_{r}(\xi_{r})=\Max{t_{r}}$.
Consider: (a) $J_r''\subset\O_{W_r}$, (b) the basic object
$(W_r,(J''_r,b''),E''_r)$ (locally at
$\xi_{r}$), corresponding, as in 5.10, to
$(W_{r},(J_{r},b),E_{r})$, so that
$\Sing(J_{r}'',b)=\MaxBar{t_{r}}$, (c)
the basic object
$(W_{r}^{(t)},(J_{r}^{(t)},b),E_{r}^{(t)})$, (d) the basic object
$(W_{r}^{(t)},((J_{r}^{(t)})'',b''),(E^{(t)}_{r})'')$ (see 5.10).
\smallskip
Then, locally at $\xi_{r}$, we have $(J_r'')^{(t)}=(J_r^{(t)})''$, where
$(J_r'')^{(t)}=J''_{r}\O_{W_{r}^{(t)}}$.

{\it Proof:}
Consider first the ideal $J_r'\subset\O_{W_r}$ which describes the set of
points where the function $\word_r$ is Maximum (see 5.10).
One can check from (3) in 6.7 and the construction of this ideal that
, locally at $\xi_{r}$, $(J_r')^{(t)}=(J_r^{(t)})'$. In fact
$\bar{J}_{r}^{(t)}=\bar{J}_{r}\O_{W_{r}^{(t)}}$ and the expressions of
the monomial parts coincide.
The same argument applies to the construction in 5.10.
(Recall that the sets $E_{r}^{-}$ and $E_{r}^{+}$ are defined in
terms of the
values $\word_{i}(\xi_{i})$ and the exceptional hypersurfaces).
This shows that $(J_r'')^{(t)}=(J_r^{(t)})''$.

\proclaim
(6.11) Lemma.
Fix the setting and assumptions of lemma~6.10, so that
\smallskip
(i) $\word_r(\xi_{r})>0$
\smallskip
(ii) $t_r(\xi_{r})=t_r^{(t)}(\xi_{r})$
(and hence and $(J_r'')^{(t)}=(J_r^{(t)})''$).
\smallskip
Then, locally at $\xi_{r}$, there exists a hypersurface $Z_r$ of $W_r$
such that:
\smallskip
(1) $Z_r$ is transversal to $W^{(t)}_{r}$ at $\xi_{r}$.
\smallskip
(2) $Z_r$ is transversal to all hypersurfaces of $E''_r$ at $\xi_{r}$.
\smallskip
(3) $I(Z_r)_{\xi_{r}}\subset\Delta_{W_r}^{b''-1}(J_r'')_{\xi_{r}}$.
\smallskip
(4) Set $Z_r^{(t)}= Z_r \cap W_r^{(t)}$, then:
(4.a) $I(Z_r^{(t)})_{\xi_{r}}\subset
\Delta_{W_r^{(t)}}^{b''-1}((J_r^{(t)})'')_{\xi_{r}}$
(in $\O_{W_{r}^{(t)}}$) and
(4.b) Locally at $\xi_{r}$,
$(E_{Z_{r}^{(t)}}^{(t)})''=(E''_{Z_{r}})^{(t)}$.
\smallskip
(5) If $\xi_{r}\in R(1)(\MaxBar{t_{r}})$ then
$\xi_{r}\in R(1)(\MaxBar{t_{r}^{(t)}})$ (where the last codimension is
considered in $W_{r}^{(t)}$)
and, locally at $\xi_{r}$,
$R(1)(\MaxBar{t_{r}})\cap W_{r}^{(t)}=R(1)(\MaxBar{t_{r}^{(t)}})$.

{\it Proof}:
Let $r_0$ be the index such that
$$\word_{r_{0}-1}(\xi_{r_{0}-1})>\word_{r_0}(\xi_{r_0})=\cdots=\word_r(\xi_r)$$
Recall, from 5.9, that the basic object $(W_r,(J_r',b'),E_r^{+})$
is the transform of the basic object
$(W_{r_0},(J_{r_0}',b'),\emptyset)$.
The basic object $(W_{r_0},(J_{r_0}',b'),\emptyset)$ satisfies
both conditions (IA)
and (LC) and we can apply lemma~6.9 to find a hypersurface $Z_{r_0}$
satisfying conditions (1), (2), (3) and (4) (formulated
now for this basic object), note
that here the set of hypersurfaces $E_{r_0}^{+}$ is empty so that (2) is
vacuous at level $ r_0 $.
Now we check that the strict transform $Z_r$ of $Z_{r_0}$ satisfies also
all four conditions for the basic object
$(W_r,(J_r',b'),E_r^{+})$ (see  5.6).

Note that, from 5.9, $b''=b'$ and $J_r'\subset J_r''$. So $Z_r$ also
satisfies condition (3) for the basic object
$(W_r,(J_r'',b''),E_{r}''=E_r^{+})$. Note also that
$ Sing(J_r'',b'')=V(\Delta_{W_r}^{b''-1}(J_r''))=\MaxBar{t_{r}}$, in
particular
$\xi_{r} \in R(1)(\MaxBar{t_{r}})$ if and only if equality holds at (3).
Now we address (4). The inclusion (4.a) follows from 6.9 since
our hypothesis asserts that $(J_r'')^{(t)}=(J_r^{(t)})''$. Equality
(4.b) follows from the definition of $E_{r}''=E_r^{+}$ in 5.9.

We finally address (5).
We assume here that $t_r(\xi_{r})=t_r^{(t)}(\xi_{r})$,
so let $ \Delta_{W_r}^{b''-1}(J_r'') $ be locally principal at $ \xi_r $
(let equality hold at (3)).
On the one hand we know from (4.a) that $I(Z_r^{(t)})_{\xi_{r}}\subset
\Delta_{W_r^{(t)}}^{b''-1}((J_r^{(t)})'')_{\xi_{r}}$. Note that,
locally $
\xi_{r} $, $I(Z_r^{(t)})=I(Z)\O_{W^{(t)}_{r}}$, so if equality
holds at (3)
then $\Delta_{W_r}^{b''-1}((J_r)'')\O_{W^{(t)}}
\subset \Delta_{W_r^{(t)}}^{b''-1}((J_r^{(t)})'')$ locally at
${\xi_{r}}$.
On the other hand the converse
inclusion always holds (see proof of 6.9). Now (5) follows.

\proclaim
(6.12) Lemma.
Fix a basic object  $(W,(J,b),E)$, a smooth pure dimensional subvariety
$W^{(t)}$ transversal to all hypersurfaces of $E$, and a point $ \xi \in
Sing(J,b)\cap W^{(t)} $. Assume that $\Max{\nu_J}=b>0 $.
Let $Z\subset W$ be a smooth hypersurface, and assume: $Z$ is
transversal to all the hypersurfaces of $E$,
$I(Z)\subset\Delta_W^{b-1}(J)$ locally at $ \xi $, and $ Z $ and
$ W^{(t)} $ are transversal at $ \xi $. Set $Z^{(t)}=Z\cap W^{(t)}$.
Assume also that $\xi\not\in R(1)(\Sing(J,b))$ and
$\xi\not\in R(1)(\Sing(J^{(t)},b))$.
Consider (as in 5.7)  the basic objects
$(Z,(C(J),b!),E_Z)$
and
$(Z^{(t)},(C(J^{(t)},b!),E^{(t)}_Z)$,
defined near $\xi$ in $W$, corresponding to
$(W,(J,b),E)$
and
$(W^{(t)},(J^{(t)},b),E^{(t)})$
respectively, as well as the
$Z^{(t)}$-ideals
$C(J)^{(t)}=C(J)\O_{Z^{(t)}}$
  and  $C(J^{(t)})$.
\smallskip
Then, $C(J^{(t)})=C(J)^{(t)}$.

{\it Proof:}
The coefficient ideal $C(J)\subset\O_{Z}$ is introduced by taking
an etale
retraction $W\longrightarrow Z$, then (locally) the ideal is
defined in terms of the
coefficients of generators
of $J$ (see part 3 in 13.1 (page 212) in [11], or
lemma~5.4 in [16]).
In our setting we consider two compatible \'etale retractions:
$W\longrightarrow Z$ and, by restriction, one
$W^{(t)}\longrightarrow Z^{(t)}$.
If $f\in\O_{W,\xi}$ let $a_{i}\in\O_{Z}$ be a coefficient of $f$.
The required
equality follows essentially from the fact that the class
$\bar{a}_{i}\in\O_{Z^{^(t)},\xi}$ is a coefficient of the class
$\bar{f}\in\O_{W^{(t)},\xi}$.
\medskip

{\bf (6.13)  Remark.} With the notation and assumptions of 6.10 and 6.11,
if $\xi_{r}\not\in
R(1)(\MaxBar{t_{r}})$, we may choose $Z_r$ and $Z_r^{(t)}$ as in 6.11
and define, according to lemma~6.12, coefficient ideals in $Z_r$ and
$Z_r^{(t)}$, and basic objects:
$$(Z_r,(C(J_r''),b''!),E''_{Z_r}) \qquad
(Z_r^{(t)},(C((J_r^{(t)})''),b''!),(E_{Z_r^{(t)}}^{(t)})'')$$

But
$(J''_{r})^{(t)}=(J_{r}^{(t)})''$ (use 6.10), and we can check that if
$\xi_{r}\not\in R(1)(\MaxBar{t_{r}})$ and $\xi_{r}\not\in
R(1)(\MaxBar{t^{(t)}_{r}})$, then

$$C(J_{r}'')^{(t)}=C((J''_{r})^{(t)})=C((J_{r}^{(t)})'')$$.

\proclaim
(6.14) Proposition.
In the setting of 6.4,
(so that $W_{0}^{(t)}\subset W_{0}$ is a smooth
pure dimensional subscheme), consider the sequence (6.1.1),
a point  $\xi_{r} \in W_{r}$,
and if $\xi_{i}\in W_{i}$ denotes the image of $\xi_{r}$
in $W_{i}$, assume that the sequence of points
$\xi_{0},\xi_{1},\ldots,\xi_{r}$
satisfies (1), (2) of 6.4 and (3) of 6.5. Assume
$\word_r(\xi_{r})>0$.  Then
$g_r(\xi_{r})\leq g_r^{(t)}(\xi_{r})$.

{\it Proof:}
By 5.11 cases (1) and (3), the resolution functions can be
expressed in the form
$g_r(\xi_{r})=(t_r(\xi_{r}),g'_r(\xi_{r}))$ and
$g_r^{(t)}(\xi_{r})=(t_r^{(t)}(\xi_{r}),{g'}_r^{(t)}(\xi_{r}))$
respectively,
where $g_r$ and $g_r^{(t)}$ are, respectively, either $\infty$ or
the resolution functions of basic objects
in one dimension less.
In fact, the functions $t$ are now as follows:
$t_r=(\word_r,n_r)$ and $t_r^{(t)}=(\word_r^{(t)},n_r^{(t)})$; where
$b . \word_r$ is the order of the ideal $\bar{J}_r$ and
$b\cdot\word_r^{(t)}$ is the order of the ideal $\bar{J}_r^{(t)}$
(see 5.6 and (3) in 6.7).
Lemma~6.8(a) together with (3) in 6.7, assert that
$\word_r(\xi_{r})\leq\word_r^{(t)}(\xi_{r})$.

If $\word_r(\xi_{r})<\word_r^{(t)}(\xi_{r})$ the result is proved.
So, assume that $\word_r(\xi_{r})=\word_r^{(t)}(\xi_{r})$.  Then, 6.9
also asserts that $n_r(\xi_{r})=n_r^{(t)}(\xi_{r})$, since, locally at
$\xi_{r}$, the set
$(E_{r}^{(t)})^{-}$ is the intersection of the hypersurfaces
$E^{-}_{r}$ with $W^{(t)}_{r}$.
We assume here that $\word_r(\xi_{r})>0$ and that
$t_r(\xi_{r})=t_r^{(t)}(\xi_{r})$.
If $ \xi_{r} \in R(1)(\MaxBar{t_{r}})$
(i.e. if $g_r(\xi_{r})=(t_r(\xi_{r}),\infty)$) then
the statement follows from (5) in Lemma 6.11, in fact
$g_r(\xi_{r})= g_r^{(t)}(\xi_{r})$ in this case.

If $\xi_{r}\not\in R(1)(\MaxBar{t_{r}}) $ and
$\xi_{r}\in R(1)(\MaxBar{t^{(t)}_{r}}) $ then
$g_r(\xi_{r})< g_r^{(t)}(\xi_{r})$.

If $\xi_{r}\not\in R(1)(\MaxBar{t_{r}}) $ and
$\xi_{r} \not\in R(1)(\MaxBar{t^{(t)}_{r}}) $ then
we choose a hypersurface $Z_{r}$ as in 6.11,
so the function $g'_{r}$ is the resolution function of the basic object
$(Z_{r},(C(J''_{r}),b''!),E''_{Z_{r}})$, and ${g'}_{r}^{(t)}$ is the
resolution function of the basic object
$(Z_{r}^{(t)},(C((J_{r}^{(t)})''),b''!),(E^{(t)}_{Z_{r}})'')$, each
corresponding to a sequence

$$(Z_{r_{0}},(C(J''_{r_{0}}),b''!),E''_{Z_{r_{0}}})
\longleftarrow\cdots\longleftarrow
(Z_{r},(C(J''_{r}),b''!),E''_{Z_{r}})$$
$$(Z^{(t)}_{r_{0}},(C((J^{(t)}_{r_{0}})''),b''!),E''_{Z^{(t)}_{r_{0}}})
\longleftarrow\cdots\longleftarrow
(Z^{(t)}_{r},(C((J^{(t)}_{r})''),b''!),E''_{Z^{(t)}_{r}})$$
for some index $r_{0}\leq r$ (see remark~6.13).
Now one can check that the two conditions of 6.6 hold for this
sequence and
points $ \xi_i $ for $r_{0}\leq i \leq r$, and hence, by induction on
the dimension of the ambient space, that $g'_r(\xi_{r})\leq
{g'}_r^{(t)}(\xi_{r})$.  So clearly $g_r(\xi_{r})\leq
g_r^{(t)}(\xi_{r})$, proving 6.14..
\smallskip

Recall, from 5.10, that in this case $g_{r}=(t_{r},g_{r})$ and
$g_{r}^{(t)}=(t_{r}^{(t)},{g'}_{r}^{(t)})$, where $g'_{r}$ is the
resolution function of
$(Z_{r},(C(J''_{r}),b''!),E''_{Z_{r}})$ and
${g'}_{r}^{(t)}$ is the one of
$(Z_{r}^{(t)},(C(J''_{r})^{(t)},b''!),(E''_{Z_{r}})^{(t)})$.
More precisely, here $g'_{r}$ is the resolution
function corresponding to the induced sequence
$$(Z_{r_{0}},(C(J''_{r_{0}}),b''!),E''_{Z_{r_{0}}})
\longleftarrow\cdots\longleftarrow
(Z_{r},(C(J''_{r}),b''!),E''_{Z_{r}})$$
for some index $r_{0}\leq r$. But for simplicity we will only say that
$g'_{r}$ is the resolution function of the last basic object.

\proclaim
(6.15) Proposition.
Fix the setting and notation as in 6.9, and assume again that
$\word_r(\xi_{r})>0$. If $g_r(\xi_{r})=g_r^{(t)}(\xi_{r})$ then
$S_{r,\alpha_{r}}$ is transversal to $W_r^{(t)}$ at the point $\xi_{r}$.

{\it Proof:}
Let $m$ be the codimension of $S_{r,\alpha_{r}}$ in $W_r$ locally at
$(\xi_{r})$.

If $m=1$ then $g_{r}(\xi_{r})=(t_{r}(\xi_{r}),\infty)$. We may
restrict ourselves
to a neighborhood of $\xi_{r}$ and assume that $\xi_{r}\in\Max{t_{r}}$.
So that $\xi_{r} \in R(1)(\MaxBar{t_{r}} ) $, in fact
$\MaxBar{g_{r}}\subset\Max{t_{r}}$.
The first assertion of 6.11(5) shows that
$\xi_{r} \in R(1)(\MaxBar{t^{(t)}_{r}} ) $.
The result follows now from (5) and (1) in 6.11.

Let $m>1$ and assume the result for $m-1$.
It follows from remark~6.13 that the resolution function of
the basic object $(Z_r^{(t)},(C(J_r'')^{(t)},b''!),(E''_{Z_r})^{(t)})$ is
the function ${g'}_r^{(t)}$.
Also the stratum $S_{r,\alpha_{r}}$ is included in $Z_{r}$, locally
at $\xi_{r}$, and $S_{r,\alpha_{r}}$ is a stratum of the
stratification defined by $g'_{r}$ in $Z_{r}$.
We argue again as in 6.14, so by induction on the dimension of the
ambient space we know that
$S_{r,\alpha_{r}}$ is transversal to
$Z_r^{(t)}$ at $\xi_{r}$. This implies the transversality of
$S_{r,\alpha_{r}}$ with $W_r^{(t)}$.

\proclaim
(6.16) Proposition.
Fix the setting and notation as in 6.14, and assume again that
$\word_r(\xi_{r})>0$.
If $S_{r,\alpha_{r}}$ is transversal to $W_{r}^{(t)}$ at $\xi_{r}$,
then $g_{r}(\xi_{r})=g^{(t)}_{r}(\xi_{r})$.

{\it Proof:}
We will proceed by induction on the local codimension $m$ of
$S_{r,\alpha_{r}}$
at the point $\xi_{r}$.

If $m=1$ then then $g_{r}(\xi_{r})=(t_{r}(\xi_{r}),\infty)$.
In fact we may consider a suitable neighborhood of $\xi_{r}$ such that
$t_{r}(\xi_{r})=\Max{t_{r}}$ and we conclude that
$\MaxBar{t_{r}}$ has a component of codimension one at $\xi_{r}$
(i.e. that  $\xi_{r} \in R(1)(\MaxBar{t_{r}} ) $) and
$S_{r,\alpha_{r}}=R(1)(\MaxBar{t_{r}}) $ locally at $\xi_{r}$.

Transversality asserts now that there exists a regular system of
parameters $x_{1},\ldots,x_{n}$ at $\O_{W_{r},\xi_{r}}$ such
that $(\bar{J}_{r})_{\xi_{r}}=(x_{1}^{b'})$ where
$b'=b\cdot\word_{r}(\xi_{r})$ and
$I(W^{(t)}_{r})=(x_{s},\ldots,x_{n})$, with $s>1$.
By (3) in 6.7 we know that the ideal $\bar{J}^{(t)}_{r}$
at $\xi_{r}$ is $(\bar{x}_{1}^{b'})$. It follows now that
$g_{r}(\xi_{r})=g^{(t)}_{r}(\xi_{r})=
\left(\left(\frac{b'}{b},p\right),\infty\right)$.

Suppose that $m>1$ and assume the result for $m-1$, note that the
order of
$\bar{J}_{r}$ is constant along $C= S_{r,\alpha_{r}}$ in a
neighborhood of $\xi_{r}$. It follows from (3) in 6.7 and
lemma~6.8(b) that
the orders of the ideals
$\bar{J}_{r}$ and $\bar{J}^{(t)}_{r}$ at $\xi_{r}$
are equal, so that the
functions $\word_{r}$ and $\word^{(t)}_{r}$ coincide at
  $\xi_{r}$.
As $W^{(t)}_{r}$ is transversal to all hypersurfaces of $E_{r}$, the
functions $n_{r}$ and $n^{(t)}_{r}$ also coincide (by (1) in 6.6);
recall here that the set $E_{r}^{-}$ is described in terms of the
values $\word_{i}(\xi_{i})$ and $E_{r}$.
So that the functions $t_{r}$ and $t^{(t)}_{r}$ take the same value
at the point $\xi_{r}$.

Recall remark~6.8.1 and that $g_{r}=(t_{r},g'_{r})$ and
$g^{(t)}_{r}=(t^{(t)}_{r},g'^{(t)}_{r})$, where $g'_{r}$ and
$g'^{(t)}_{r}$ are the resolution functions of basic objects in smaller
dimensions ($Z_{r}$ and $Z_{r}^{(t)}$ respectively).

Note also that $S_{r,\alpha_{r}}\subset Z_{r}$ locally at $\xi_{r}$
and $S_{r,\alpha_{r}}$ is also one stratum of the stratification of
$g'_{r}$ in $Z_{r}$.
By hypothesis we conclude that $Z_r^{(t)}$ is transversal to
$S_{r,\alpha_{r}}$,
so by induction ${g'}_r(\xi_{r})={g'}_r^{(t)}(\xi_{r})$.
It follows that $g_r(\xi_{r})=g_r^{(t)}(\xi_{r})$.

\proclaim
(6.17) Proposition.
Fix the setting and notation as in 6.14, and assume now
that $\word_r(\xi_{r})=0$. Then, $g_{r}(\xi_{r})=g_{r}^{(t)}(\xi_{r})$
and $S_{r,\alpha_{r}}$ is
transversal to $W_{r}^{(t)}$ at $\xi_{r}$.

{\it Proof:}
Note that it follows from 5.11 case (2), that
$g_{r}(\xi_{r})=\Gamma_{r}(\xi_{r})$ and
$g_{r}^{(t)}(\xi_{r})=\Gamma_{r}^{(t)}(\xi_{r})$.
But $\Gamma_{r}$ and $\Gamma_{r}^{(t)}$ depend only on the expression
of the monomial parts of $J_{r}$ and $J_{r}^{(t)}$, respectively,
and we know by condition (3) in 6.7 that these expressions are the same.

The set $S_{r,\alpha_{r}}$ is, locally, some
intersection of the hypersurfaces of $E_{r}$, so that it is
transversal to $W_{r}^{(t)}$, since $W_{r}^{(t)}$ is transversal to
all hypersurfaces of $E_{r}$.
\medskip

{\bf (6.18)} {\it Proof of Theorem 6.4.}
It is by induction on $r$ (the order of compatibility), and will
break down
into different cases.
Note that for $r=0$, then conditions (1)
and (2) of 6.4 are void, in fact they are formulated for indices
$i=0,\ldots,r-1$. So, for $r=0$,
our local assumptions will hold, and Theorem~6.4 follows
from propositions~6.14, 6.15 and 6.16.
\medskip

Assume that $r>0$ and that Theorem~6.4 holds for $r'<r$.
We claim that we may reduce the proof of our theorem to the case
where $\xi_{r-1}\in C_{r-1}$.
In fact, if
$\xi_{r-1}$ does not belong to the center $C_{r-1}$ then by property
p0 of 5.3 we may identify $\xi_{r}$ and $\xi_{r-1}$ and the conditions
follow from Theorem~6.4, applied for $r-1$. So set
$\xi_{0},\xi_{1},\ldots,\xi_{r}$ as before, assume now that
$\xi_{r-1}\in C_{r-1}$ and consider within the set of indices ${0 \leq j
\leq r-1}$ the subset $ F $ of those indeces $j_0 $ such that
$\xi_{j_0}\in
C_{j_0}$ (so
$ r-1 \in F $). Let us check now that if conditions (C1) and (C2) of 6.1
hold, then also local conditions (1) and (2) of 6.6 hold. Fix an
index $ k
$, so that ${0 \leq k \leq r-1}$. If $\xi_{k}\in C_{k}$, then (2)
of 6.4 (for
index $k$) follows from (C1) and (C2) of 6.1; and (1) of 6.6 (for
index $k$)
follows by our induction hypothesis. Now let $j $ be an arbitrary index,
${0\leq j \leq r-1}$. Since $\xi_{r-1}\in C_{r-1}$, one can check that
there is a well
defined $j_0 $, so that $j \leq j_0 $,  $\xi_{j_{0}}\in C_{j_{0}}$,
and the
transformations define a local isomorphism of
$\xi_{j}\in W_j $ (of $\xi_{j}\in W^{(t)}_j $) with
$\xi_{j_0}\in W_{j_0} $ (of $\xi_{j_0} \in W^{(t)}_{j_0} $). Hence,
by p0 of 5.3, both
conditions (1) and (2) of 6.6 hold.
\medskip

Now Theorem~6.4 will follow, for $r>0 $, from propositions~6.14,
6.15, 6.16
and 6.17.
\vskip 1cm

{\bf REFERENCES.}
\medskip

\parindent=14mm
\frenchspacing

\item{[1]} S.~S.~Abhyankar. {\it A criterion of equisingularity},  
Amer. J.
Math  90 (1968), 342-345.
\smallskip

\item{[2]} S.~S.~Abhyankar. {\it Remarks on Equisingularity}, Amer.  
J. Math 90
(1968), pp.108-144.
\smallskip

\item{[3]} D. Abramovich and A.J. de Jong.
{\it Smoothness, semistability and toroidal geometry}, Journal of
Algebraic
Geometry, 6 (1997), 789-801.
\smallskip

\item{[4]} D. Abramovich and J. Wang, {\it Equivariant
resolution of singularities in characteristic 0},
Mathematical Research Letters, 4
(1997), 427-433.
\smallskip

\item{[5]} P. Berthelot, {\it Alt\'erations de variet\'es
alg\'ebriques}, S\'eminaire Bourbaki, Expos\'e No. 815,
Ast\'erisque, 241 (1997),  273-311.

\item{[6]} E. Bierstone and P. Milman, {\it Canonical
desingularization
  in characteristic zero by blowing-up the maximal
strata of a local invariant}, Inv. Math. 128 (1997), 207-302.
\smallskip

\item{[7]} F. Bogomolov and T. Pantev, {\it Weak Hironaka
Theorem}, Mathematical Research Letters 3 (1996),  299-307.

\item{[8]}
G.~Bodn\'{a}r and J.~Schicho.
{\it A Computer Program for the Resolution of Singularities}
in {\it Resolution of Singularities. A research textbook in tribute to
Oscar Zariski}, (H. Hauser, J. Lipman, F.
Oort and A. Quir\'os editors)
Progress in Mathematics, Volume 181, Birkh\"auser Verlag,
Basel, 2000.
\smallskip

\item{[9]}
G.~Bodn\'{a}r and J.~Schicho.
{\it Automated Resolution of Singularities for Hypersurfaces}
in {\it Applications of Groebner Bases}.
To appear. Available at:

{\tt http://www.risc.uni-linz.ac.at/projects/basic/adjoints/blowup}
\smallskip

\item{[10]} S. Encinas and O. Villamayor, {\it Good points and
constructive resolution of singularities}, Acta Math. Vol. 181
(1998), pp. 109--158.
\smallskip

\item{[11]}
S.~Encinas and O.~Villamayor.
{\it A Course on Constructive Desingularization and Equivariance}
in {\it Resolution of Singularities. A research textbook in tribute to
Oscar Zariski}, (H. Hauser, J. Lipman, F.
Oort and A. Quir\'os editors)
Progress in Mathematics, Volume 181, Birkh\"auser Verlag,
Basel, 2000.
\smallskip

\item{[12]}
S.~Encinas and O.~Villamayor.
{\it A new theorem of desingularization over fields of characteristic
zero}.
Preprint.
\smallskip

\item{[13]}
T.~Gaffney and D.~Massey.
{\it Trends in equisingularity theory}, in {\it Singularity theory},
Bill Bruce and David Mond, editors.
London Math. Soc. lecture series 263, 207-248, Cambridge
Univ. Press, Cambridge.

\item{[14]} R. Hartshorne, {\it Algebraic Geometry},
Graduate Texts in Mathematics 52,
Springer Verlag, New York (1983).
\smallskip

\item{[15]} H. Hironaka, {\it Resolution of singularities of an
algebraic variety over a field of characteristic zero I-II}, Ann. Math.,
79 (1964), 109-326.
\smallskip

\item{[16]} H. Hironaka, {\it Idealistic exponent of a
singularity}, Algebraic Geometry. The John Hopkins centennial lectures,
Baltimore, Johns Hopkins University Press (1977), pp. 52--125.
\smallskip

\item{[17]} A.J. de Jong, {\it Smoothness, semi-stability and
alterations}, Pub. Math. I.H.E.S. no. 83 (1996), 51-93.
\smallskip

\item{[18]} S.~Kleiman. {\it Equisingularity, Multiplicity and
Dependence}. Preprint 1998.\break alg-geom/9805062.

\item{[19]} H.~Laufer.
{\it Strong simultaneous resolution for surface singularities}, in
{\it Complex Analytic Singularities}, North-Holland,
Amsterdam-New York (1987), 207-214.
\smallskip

\item{[20]} J.~Lipman.
{\it Equisingularity and Simultaneous Resolution of Singularities}
in {\it Resolution of Singularities. A research textbook in tribute to
Oscar Zariski}, (H. Hauser, J. Lipman, F.
Oort and A. Quir\'os editors)
Progress in Mathematics, Volume 181, Birkh\"auser Verlag,
Basel, 2000.
\smallskip

\item{[21]} H. Matsumura,
{\it Commutative Algebra (second Edition)}
Mathematics Lecture Note Series,
Benjamin, Cumming Publishing Company (1980).
\smallskip

\item{[22]} D.~Mumford.
{\it Geometric Invariant Theory}.
Ergebnisse der Mathematik, Vol. 34, Springer-Verlag, Berlin (1965).
\smallskip

\item{[23]} A. Nobile, {\it Equisingular stratifications of
parameter spaces of families of embedded curves}, to appear in {\it
Proc. Conference in honor of O. Villamayor, Sr.},  Acad. Nac.
Ciencias Cordoba, publicacion No 65.
\smallskip

\item{[24]} A. Nobile and O.E. Villamayor, {\it Equisingular
stratifications associated to families of planar ideals}, J. Algebra
193 (1997), 239-259 .
\smallskip

\item{[25]} A. Nobile and O.E. Villamayor, {\it Arithmetic
families of smooth surfaces and equisingularity of embeded schemes},
Manuscripta Math. 100, 173-196 (1999), 173-196 .
\smallskip

\item{[26]} S.~A.~Stromme.
{\it Elementary introduction to representable functors and
Hilbert Schemes}.
Parameter Spaces. Banach Center Publications, Vol 36.
Institute of Mathematics. Polish Academy of Science
Warszawa 1996
\smallskip

\item{[27]} B.~Teissier.
{\it R\'esolution simultan\'ee},
in {\it S\'eminaire sur les Singularit\'es de Surfaces}, (M.
Demazure, H. Pinkham and B. Teissier editors)
Lecture Notes in Mathematics,Vol. 777, Springer Verlag, Berlin (1980),
82-146.
\smallskip

\item{[28]} O. E. Villamayor, {\it Constructiveness of Hironaka's
resolution}, Ann. Scient. Ec. Norm. Sup., $4^{\rm{e}}$ s\'erie, 22
(1989), 1-32.
\smallskip

\item{[29]} O. E. Villamayor, {\it On equiresolution and a question
of Zariski}, preprint.
\smallskip

\item{[30]} J.~Wahl.
{\it Equisingular deformations of normal surface singularities},
Ann. Math. 104 (1976), 325-356.
\smallskip

\item{[31]} O.~Zariski.
{\it Some open questions in the theory of
singularities}, Bull. Amer. Math. Soc. 77 (1971),481-491.
\smallskip

\item{[32]} O.~Zariski.
{\it Foundations of a general theory of equisingularity
  on $r$-dimensional algebroid and algebraic varieties,
of embedded dimension $r+1$}, Amer. J. Math. 101 (1979), 453-514.
(reprinted in [33], pp. 573-634; and summarized in ibid, pp. 635-651.)
\smallskip

\item{[33]} O.~Zariski. Collected Papers, vol. IV,
MIT Press, Cambridge, Mass., 1979.

\bigskip
\bigskip

{S. Encinas

\nobreak
Dep. de Matem\'atica Aplicada Fundamental.

\nobreak
E.T.S. Arquitectura. Universidad de Valladolid.

\nobreak
Avda. Salamanca s/n

\nobreak
47014 Valladolid.

\nobreak
Spain.

\nobreak
sencinas@maf.uva.es}
\medskip

{A. Nobile

\nobreak
Dep. of Mathematics.

\nobreak
Louisiana State University.

\nobreak
Bat\^on Rouge.

\nobreak
Louisiana 70803

\nobreak
USA.

\nobreak
nobile@marais.math.lsu.edu}
\medskip

{O. Villamayor

\nobreak
Dep. de Matem\'aticas.

\nobreak
Universidad Aut\'onoma de Madrid.

\nobreak
Ciudad Universitaria de Canto Blanco.

\nobreak
28049 Madrid.

\nobreak
Spain.

\nobreak
villamayor@uam.es}

\bigskip

\bye